\documentclass[11pt]{amsart}
\usepackage{amsfonts}
\usepackage{amsmath}
\usepackage{amssymb} 
\usepackage{amsthm}
\usepackage{ color, ulem}
\usepackage{tikz-cd}
\usepackage{float}
\input xy
\xyoption {all}

\newcommand{\C}{\mathbb{C}}

\newcommand{\Quot}{\mathsf{Quot}}
\newcommand{\Hilb}{\mathsf{Hilb}}
\newcommand{\Obs}{\textnormal{Obs}}
\newcommand{\vir}{\textnormal{vir}}

\renewcommand{\dim}{\textnormal{dim}}

\newcommand{\comment}[1]{}
\newtheorem{theorem}{Theorem}
\newtheorem {lemma}[theorem]{Lemma}
\newtheorem{conjecture}[theorem]{Conjecture}
\newtheorem{question}[theorem]{Question}

\newtheorem {proposition}[theorem]{Proposition}
\theoremstyle{definition}
\newtheorem {example}[theorem]{Example} 
\theoremstyle {definition} 
\newtheorem{remark}[theorem]{Remark}

\begin{document}

\baselineskip=17.5pt

\title [The virtual $K$-theory of Quot schemes of surfaces]
{The virtual $K$-theory of Quot schemes of surfaces}

\author[Arbesfeld]{Noah Arbesfeld}
\address{Imperial College London, Department of Mathematics, London, UK}
\email{n.arbesfeld@imperial.ac.uk }
\author[Johnson]{Drew Johnson}
\address{Department of Mathematics, ETH Z\"urich}
\email {jared.johnson@math.ethz.ch}
\email {d.johnson@mathematics.byu.edu}
\author[Lim]{Woonam Lim}
\address{Department of Mathematics, University of California, San Diego}
\email{w9lim@ucsd.edu}
\author[Oprea]{Dragos Oprea}
\address{Department of Mathematics, University of California, San Diego}
\email {doprea@math.ucsd.edu}
\author[Pandharipande]{Rahul Pandharipande}
\address{Department of Mathematics, ETH Z\"urich}
\email {rahul@math.ethz.ch}

\begin{abstract}

We study virtual invariants of Quot schemes parametrizing quotients of dimension at most $1$ of the trivial sheaf of rank $N$ on nonsingular projective surfaces. We conjecture that the generating series of virtual $K$-theoretic invariants are given by rational functions. We prove rationality for several geometries including punctual quotients for all surfaces and dimension $1$ quotients for surfaces $X$ with $p_g>0$. We also show that the generating series of virtual cobordism classes can be irrational. 

Given a $K$-theory class on $X$ of rank $r$, we associate natural series of
virtual Segre and Verlinde numbers. We show that the
Segre and Verlinde series match in the following cases:
  \begin{itemize}\item [(i)] Quot schemes of dimension 0 quotients, \item [(ii)] Hilbert schemes of points and curves over surfaces with $p_g>0$, \item [(iii)] Quot schemes of minimal elliptic surfaces for quotients supported on fiber classes. \end{itemize}  Moreover, for punctual quotients of the trivial sheaf of rank $N$, we prove a new symmetry of the Segre/Verlinde series exchanging  $r$ and $N$. The Segre/Verlinde statements have analogues for punctual Quot schemes over curves. 
\end{abstract}
\date{February, 2021}
\maketitle 

\setcounter{tocdepth}{1}
\tableofcontents

\section{Introduction} 

\subsection{Overview}\label{y1}
While moduli spaces of sheaves on higher dimensional varieties
rarely carry $2$-term perfect obstructions theories, 
moduli spaces of sheaves on varieties of dimension at most 3
often have well-defined virtual fundamental classes. In many cases, the resulting virtual invariants have a rich structure, reflecting the underlying geometry of the moduli spaces used in their definition. We refer the reader to \cite {PT3} for an introduction to sheaf counting methods in enumerative geometry. 

An example whose virtual geometry can be studied effectively is the Quot scheme
of $1$-dimensional quotients over surfaces. If $X$ is a nonsingular projective 
surface and $\beta\in H_2(X,\mathbb{Z})$ is an effective curve class,
we let  $\mathsf {Quot}_X(\mathbb C^N, \beta, n)$ parametrize
short exact sequences \begin{equation}\label{qqqq}
  0\to S\to \mathbb C^N\otimes \mathcal O_X\to Q\to 0\, ,
  \end{equation}
where $$\text{rank }Q=0\, , \,\,\, c_1(Q)=\beta\, ,\,\,\, \chi(Q)=n\, .$$ 
By \cite {MOP}, $\mathsf {Quot}_X(\mathbb C^N, \beta, n)$ carries a canonical
$2$-term perfect obstruction theory and a virtual fundamental class of
dimension $$\text{vdim}=\chi(S, Q)=Nn+\beta^2\, .$$ The virtual fundamental class of Quot schemes over curves had been constructed previously in \cite{mo1}.
Due to connections with Seiberg-Witten theory,
the virtual fundamental class of Quot schemes over surfaces with $p_g>0$ is
typically more accessible, while the case $p_g=0$ is less understood. 

For the Quot schemes of $1$-dimensional quotients over surfaces,
the generating series of the following invariants have been
considered in \cite {OP, L, JOP} respectively: \begin{itemize}
\item [(i)] virtual Euler characteristics,
\item [(ii)] virtual $\chi_y$-genera,
\item [(iii)] descendent invariants.
\end{itemize} 
Conjecturally, series (i) -- (iii) are always given by rational functions.
Rationality was shown for arbitrary surfaces when $\beta=0$, see \cite{JOP, L, OP}. 
For non-zero curve classes, rationality was confirmed for (i) and (ii) for all surfaces with $p_g>0$, see \cite {L, OP}.
Furthermore, series (i) and (iii) were also proven to be rational for simply connected surfaces with $p_g=0$ when $N=1$ in \cite{JOP}. 

We refine here the techniques of \cite {JOP, L, OP}
to study the virtual $K$-theory of the Quot scheme $\mathsf{Quot}_{X}(\mathbb C^N, \beta, n)$.
A central result is the rationality of natural series of
virtual $K$-theoretic invariants for many geometries. The methods also
allow us
to prove several new identities and symmetries satisfied by the virtual Segre and Verlinde series. Along the way, we take up the series of cobordism classes, and study the virtual structure sheaf for punctual quotients. 

\subsection{Seiberg-Witten invariants}
Let $X$ be a nonsingular projective surface.
Curve classes $\beta$ of Seiberg-Witten length $N$  were introduced in \cite {L}.
By definition,
$\beta\in H_2(X,\mathbb{Z})$ is of Seiberg-Witten length $N$ if 
for all effective decompositions $$\beta=\beta_1+\ldots+\beta_N$$ such that $\mathsf{SW}(\beta_i)\neq 0\text{ for all }i$, we have $$\beta_i\cdot (K_X-\beta_i)=0 \text { for all }i\, .$$
Here, $\mathsf{SW}(\beta_i)\in H_{\star}(\text{Pic}(X))=\wedge^{\star} H^1(X)$ denotes the Seiberg-Witten invariant of \cite {DKO} constructed via the Hilbert scheme of curves $\mathsf {Hilb}_{\beta_i}(X)$ of class $\beta_i$ and the associated
Abel-Jacobi map $$\mathsf {AJ}:\mathsf {Hilb}_{\beta_i}(X)\to \text{Pic}^{\beta_i}(X).$$
When the condition $\beta_i\cdot (K_X-\beta_i)=0$ is satisfied, the Hilbert scheme of curves admits a virtual fundamental class of dimension zero whose length is $$\mathsf {SW}(\beta_i)=\text{deg}\, \left[\mathsf {Hilb}_{\beta_i}(X)\right]^{\mathrm{vir}}\in \mathbb{Z}\, .$$ As noted in \cite {DKO, L}, examples of curve classes of Seiberg-Witten length $N$, for all $N$, include:
\begin{itemize}
\item [(i)] $\beta=0$ for all surfaces,
\item [(ii)] arbitrary curve classes $\beta$ for surfaces with $p_g>0$, 
\item [(iii)] curve classes supported on fibers for relatively minimal elliptic surfaces.
\end{itemize}

\subsection{K-theory}\label{kthry} 
For a scheme $S$ with a $2$-term perfect obstruction theory, let $\mathcal O_S^{\vir}$ denote the virtual structure sheaf, see \cite {BF, CFK, Lee}. Given a $K$-theory class $V\to S$,
write $$\chi^{\text{vir}}(S, V)=\chi(S, V\otimes \mathcal O_S^{\text{vir}})\, .$$ 

For $\alpha \in K^0(X)$, we define the tautological classes 
$$\alpha^{[n]}= {\textnormal{R}}\pi_{1*}(\mathcal{Q} \otimes \pi_2^*\alpha) \in K^0(\mathsf {Quot}_X(\mathbb C^N, \beta, n))\, .$$
Here, $$\mathcal{Q}\to\mathsf {Quot}_X(\mathbb C^N, \beta, n) \times X$$ is the universal quotient, and $\pi_1, \pi_2$ are the two projections.  

The Riemann-Roch numbers of exterior and symmetric powers of tautological sheaves over the Hilbert schemes of points $X^{[n]}$ (with its usual geometry) were studied in \cite {A, D, EGL, K, Sc, Z}, among others. Closed-form expressions are difficult to write down. By working with the virtual class, we cover more general Quot scheme geometries and obtain
 answers that satisfy simpler structural results. We first prove

 \begin{theorem} \label{n2} Fix integers $k_1, \ldots, k_{\ell}\geq 0$, and let $\alpha_1, \ldots, \alpha_{\ell}$ be $K$-theory classes on a nonsingular projective surface $X$.
   If $\beta\in H_2(X,\mathbb{Z})$ is a curve class of Seiberg-Witten length $N$, the series $$\mathsf Z_{X, N, \beta}^{K} (\alpha_1, \ldots, \alpha_{\ell}\,|\,k_1, \ldots, k_{\ell})=
   \sum_{n\in \mathbb{Z}} q^n \chi^{\textnormal{vir}}\left(\mathsf{Quot}_{X}(\mathbb C^N, \beta, n), \wedge^{k_1} \alpha_1^{[n]}\otimes \ldots \otimes \wedge^{k_{\ell}} \alpha_{\ell}^{[n]}\right)$$ is the Laurent expansion of a rational function in $q$. 
\end{theorem}
 
By taking $\alpha$ of negative rank, we also access the symmetric powers of the tautological sheaves via the identities $${\wedge}_{\,t}\alpha=\sum_{k} t^k \wedge^k \alpha\, ,\ \ \, \mathsf{S}_t (\alpha) = \sum_{k} t^k\, \text{Sym}^k \alpha\, , \,\,\,\,{\wedge}_{-t}(-\alpha)= \mathsf{S}_t (\alpha)\, .$$ 

Theorem \ref{n2} applies to all curve classes on surfaces with $p_g>0$. The case of surfaces with $p_g=0$ is more complicated, and our results are not as general. Nonetheless, when $N=1$, we show 

\begin{theorem}\label{n3} For all nonsingular simply connected surfaces $X$ with $p_g=0$ and all  $\alpha\in K^0(X)$, the series $$\mathsf Z_{X, 1, \beta}^{K} =
   \sum_{n\in \mathbb{Z}} q^n \chi^{\textnormal{vir}}\left(\mathsf{Quot}_{X}(\mathbb C^1, \beta, n),\, \alpha^{[n]}\right)$$ is the Laurent expansion of a rational function in $q$. 
\end{theorem}

Based on Theorems \ref{n2} and \ref{n3}, we formulate the following 

 \begin{conjecture} \label{c2} Let $k_1, \ldots, k_{\ell}\geq 0$ be integers, and let $\alpha_1, \ldots, \alpha_{\ell}$ be $K$-theory classes on a nonsingular projective surface $X$.
For every $\beta\in H_2(X,\mathbb{Z})$, the series $$\mathsf Z_{X, N, \beta}^{K} (\alpha_1, \ldots, \alpha_{\ell}\,|\,k_1, \ldots, k_{\ell})=
   \sum_{n\in \mathbb{Z}} q^n \chi^{\textnormal{vir}}\left(\mathsf{Quot}_{X}(\mathbb C^N, \beta, n), \wedge^{k_1} \alpha_1^{[n]}\otimes \ldots \otimes \wedge^{k_{\ell}} \alpha_{\ell}^{[n]}\right)$$ is the Laurent expansion of a rational function in $q$. 
\end{conjecture}

Several other variations are possible. We can dualize some of the factors of the tensor product. We could also consider twists by the virtual cotangent bundle \begin{align*}\mathsf {\widetilde Z}_{X, N, \beta}^{K} (\alpha_1, \ldots, \alpha_{\ell}\,|\,k_1, &\ldots, k_{\ell})=\\&=\sum_{n\in \mathbb{Z}} q^n \chi^{\textrm{vir}}\left(\mathsf{Quot}_{X}(\mathbb C^N, \beta, n), \wedge^{k_1} \alpha_1^{[n]}\otimes \ldots \otimes \wedge^{k_{\ell}} \alpha_{\ell}^{[n]}\otimes \wedge_y \Omega^{\textrm{vir}}\right)\, .\end{align*} These variations can be studied by the same methods, yielding rational functions in $q$, though we do not explicitly record the answers.

It is natural to ask about the locations and orders of the poles of the $K$-theoretic
generating series. 
 In case $N=1$,  we have a complete answer. 
 \begin{theorem} \label{n1}
   If $\beta\in H_2(X,\mathbb{Z})$ is a curve class of Seiberg-Witten length $N=1$, the shifted
   series $$\overline{\mathsf Z}_{X, 1, \beta}^{K} (\alpha_1, \ldots, \alpha_{\ell}\,|\,k_1, \ldots, k_{\ell})=
 q^{\beta\cdot K_X}  \sum_{n\in \mathbb{Z}} q^n \chi^{\textnormal{vir}}\left(\mathsf{Quot}_{X}(\mathbb C^1, \beta, n), \wedge^{k_1} \alpha_1^{[n]}\otimes \ldots \otimes \wedge^{k_{\ell}} \alpha_{\ell}^{[n]}\right)$$ is the Laurent expansion of a rational function 
with a pole only at $q=1$ of order at most $$2(k_1+\ldots+k_{\ell})\, .$$ 
\end{theorem}
We formulate the following non-virtual analogue of Theorem \ref{n1}.
\begin{question}\label{q5} For integers $k_1, \ldots, k_{\ell}\geq 0$ and $K$-theory classes $\alpha_1, \ldots, \alpha_{\ell}$, is the series $$\sum_{n=0}^{\infty} q^n \chi\left(X^{[n]}, \wedge^{k_1} \alpha_1^{[n]}\otimes \ldots \otimes \wedge^{k_{\ell}} \alpha_{\ell}^{[n]}\right)$$ the Laurent expansion of a
  rational function in $q$?
\end{question} 
\noindent Explicit calculations in \cite [Section 6]{A}, \cite [Section 5]{EGL}, \cite [Section 8]{K} , \cite [Section 5]{Sc}, \cite [Section 7]{Z} answer Question \ref{q5} in the affirmative for several values of the parameters $\ell$, $k_i$ and $\text{rank }\alpha_i$. We will investigate the general case in future work.

\begin{example} \label{exa1} The simplest case of Theorem \ref{n2} occurs when $k_i=0$ for all $i$. By the results of \cite {L},  \begin{equation}\label{taut0}\mathsf Z^{K}_{X, N}=
    \sum_{n\in \mathbb{Z}} q^n \chi^{\textrm{vir}}\left(\mathsf{Quot}_{X}(\mathbb C^N, n), \mathcal O\right)=1\,
    ,\end{equation} and, in general, $$\mathsf Z^{K}_{X, N, \beta}=\sum_{n\in \mathbb{Z}} q^n \chi^{\textrm{vir}}\left(\mathsf{Quot}_{X}(\mathbb C^N, \beta, n), \mathcal O\right)$$ is rational whenever $\beta$ is of Seiberg-Witten length $N$. 
\end{example}
\begin{example}\label{exa} Let $\beta=0$, and take $\ell=1$, $k_1=1$. To illustrate Theorems \ref{n2} and \ref{n1}, we compute\footnote{Non-virtually, we can compare \eqref{taut1} with Proposition 5.6 in \cite {EGL}:
$$\sum_{n=1}^{\infty} q^{n-1} \chi  (X^{[n]},\alpha^{[n]}) = \frac{ \chi(\alpha)}{(1-q)^{\chi(\mathcal O_X)}}\, .$$}
\begin{equation}\label{taut1}
\sum_{n=0}^{\infty} q^n \chi^{\textrm{vir}} (X^{[n]}, \alpha^{[n]})=- \text{rk } \alpha \cdot K_X^2\cdot \frac{q^2}{(1-q)^2} - \langle K_X, c_1(\alpha) \rangle \cdot \frac{q}{1-q}\, .
\end{equation}
Here, $\langle,\rangle$ is the intersection pairing on $X$.
For $N>1$, we find the surprising identity
\begin{equation}\label{taut2}\sum_{n=0}^{\infty} q^n \chi^{\textrm{vir}} (\mathsf{Quot}_{X}(\mathbb C^N, n), \alpha^{[n]})= N\,\sum_{n=0}^{\infty} q^n \chi^{\textrm{vir}} (X^{[n]}, \alpha^{[n]})\, .\end{equation}
Since the above series vanish for $K$-trivial surfaces, the following question
is natural.

\begin{question} What are the corresponding series for the reduced tangent-obstruction theory of Quot schemes of $K3$ surfaces? 
\end{question} 
\end{example}

\begin{example}\label{exc} Perhaps the simplest example with $\beta\neq 0$ is the case of relatively minimal elliptic surfaces $X\to C$ with $\beta$ supported on fibers. We have 
  \begin{equation}\label{taut3}\sum_{n=0}^{\infty} q^n \chi^{\textrm{vir}} (\mathsf{Quot}_{X}(\mathbb C^N, \beta,  n),  \alpha^{[n]})=\mathsf {sw}_{\beta} \cdot \left(-N 
      \langle K_X, c_1(\alpha)\rangle \frac{q}{1-q}+\langle \beta, c_1(\alpha)\rangle \frac{1}{1-q}\right) \end{equation} where $$\mathsf {sw}_{\beta}=\sum_{\beta_1+\ldots+\beta_N=\beta} \mathsf {SW}(\beta_1)\cdots \mathsf {SW}(\beta_N)\, .$$ The Seiberg-Witten invariants of fiber classes are known by \cite[Proposition 5.8]{DKO}.
\end{example}

\begin{example} \label{exd} For a more complicated example, let $X$ be a minimal surface of general type with $p_g>0$, and let $\beta=K_X$. We find \begin{equation}\label{taut4}\sum_{n\in \mathbb Z} q^n \chi^{\textrm{vir}} (\mathsf{Quot}_{X}(\mathbb C^N, K_X,  n), \alpha^{[n]})=\mathsf {SW}(K_X)\cdot \left(-\frac{N}{1-q}\right)^{K_X^2}\left(\sum_{i=1}^{N} z_i^{-K_X^2}\cdot \mathsf{p}_{\alpha}(z_i)\right).\end{equation} Here, $z_1, \ldots, z_N$ are the roots of the equation $$z^N-q(z-1)^N=0.$$ 
 In general, $\mathsf{p}_{\alpha}$ is an explicit rational function, but the simplest expressions occur when $\text{rank }\alpha=0$:
$$\mathsf{p}_{\alpha}(z)=\langle K_X, c_1(\alpha)\rangle \left(\frac{1-(N+1)q}{1-q}-z\right).$$
\end{example}

\subsection{Stable pairs} The Quot schemes studied here can be compared to the moduli space of higher rank stable pairs on surfaces $X$,
$$\mathbb C^N\otimes \mathcal O_X\stackrel{s}{\longrightarrow} F,\quad [F]=\beta\, , \quad \chi(F)=n\, ,$$ whose virtual fundamental class was constructed in \cite{Lin}. The tautological classes $\alpha^{[n]}$ can be defined analogously.
It is natural to inquire whether the above rationality results still hold.
 
Increasing the dimension, let $X$ be a smooth projective threefold, and let $P_n(X, \beta)$ denote the moduli space of stable pairs 
 $$\mathcal O_X\stackrel{s}{\longrightarrow} F,\quad [F]=\beta, \quad \chi(F)=n,$$
 as defined in \cite {PT1}. In cohomology, the series of descendant
 invariants{\footnote{The descendent series \eqref{lss559} is in
     a slightly different
     form than the descendent series in the stable pairs references, but the
     rationality of \eqref{lss559} follows from the descendent study in the
     references by an application of Grothendieck-Riemann-Roch.}}
 \begin{equation}
   \label{lss559}
 \sum_{n\in \mathbb{Z}} q^{n} \cdot \int_{\left[P_n(X, \beta)\right]^{\text{vir}}}
  \mathsf{ch}_{k_1}(\alpha_1^{[n]}) \ldots \mathsf{ch}_{k_\ell}(\alpha_\ell^{[n]})
  \,
  \end{equation}
 is conjectured to be rational \cite{MNOP, P, PP1, PP2, PT1, PT2}. Rationality is proven for toric threefolds in \cite {PP1, PP2} via localization of the virtual class \cite{GP}. The analogous conjecture in $K$-theory reads
\begin{conjecture} \label{pss45} For all integers $k_1, \ldots, k_{\ell}\geq 0$, and $K$-theory classes $\alpha_1, \ldots, \alpha_{\ell}$ on a smooth projective threefold $X$, the series 
$$\mathsf Z_{X, \beta}^{K} (\alpha_1, \ldots, \alpha_{\ell}\,|\,k_1, \ldots, k_{\ell})=
   \sum_{n\in \mathbb{Z}} q^n \chi^{\textnormal{vir}}\left(P_n(X, \beta), \wedge^{k_1} \alpha_1^{[n]}\otimes \ldots \otimes \wedge^{k_{\ell}} \alpha_{\ell}^{[n]}\right)$$ is the Laurent expansion of a rational function in $q$. 
 \end{conjecture} 
 For stable pairs on threefolds, rationality of the descendent series in cohomology {\it implies}\footnote{See also \cite{Smirnov} for a direct approach to the rationality of the
1-leg descendent vertex in $K$-theory via quasimaps.} the rationality of Conjecture \ref{pss45}  for $K$-theory via the virtual Hirzebruch-Riemann-Roch theorem \cite{FG}. Indeed, the $K$-theoretic invariants can be expressed solely in terms of the Chern characters $\text{ch}_k(\alpha^{[n]})$ appearing in the descendent series. The Chern classes of the virtual tangent bundle appearing in the calculation
 can be written in terms of descendent invariants in a form which does not
 depend on $n$, see \cite[Section 3.1]{Sh}. The result crucially
 relies on the independence of 
 the virtual dimension of $P_n(X, \beta)$ on $n$. 

 For Quot schemes of surfaces, the descendent rationality of \cite{JOP}
 does {\it not} obviously imply Conjecture \ref{c2}.
 The virtual dimension of the Quot scheme on a surface grows with $n$.
 Therefore, we cannot bound the degree of the descendent classes appearing in Hirzerbruch-Riemann-Roch independently of $n$, and hence the $K$-theoretic series are not expressed in terms of finitely many descendent series
\begin{equation}\label{desjop} \sum_{n\in \mathbb{Z}} q^{n} \cdot \int_{\left[\mathsf {Quot}_X(\mathbb C^N, \beta, n)\right]^{\text{vir}}}
  \mathsf{ch}_{k_1}(\alpha_1^{[n]}) \ldots \mathsf{ch}_{k_\ell}(\alpha_\ell^{[n]})\, c(T^{\text{vir}}\mathsf {Quot})
  \,\end{equation}
 of \cite{JOP}. Nevertheless, we will remark in Section \ref{s3} that the methods used in the proof of Theorem \ref{n2} also establish rationality of the Quot scheme descendent series for all surfaces with $p_g>0$, partially answering Conjecture $2$ of \cite{JOP}.

\subsection{Cobordism} 
It is natural to inquire whether the rationality results proven in cohomology or $K$-theory also hold at the level of virtual cobordism.{\footnote{For
    the parallel theory of stable pairs over $3$-folds, both the
    series of invariants and the
    series of cobordism classes \cite{Sh} are conjectured to be rational.}}
We show that this is not the case:
\begin{theorem}\label{p1} The series $$\mathsf Z^{\mathsf{cob}}_{X, \beta}=\sum_{n\in \mathbb Z} q^n \left[\mathsf{Quot}_{X}(\mathbb C^1, \beta, n)\right]^{\textnormal{vir}}_{\textnormal{cob}}\in \Omega^{\star}((q))$$ is not given by a rational function, in general.  
\end{theorem} 
In particular, we show that the virtual Pontryagin series is given by an explicit algebraic irrational function when $N=1$ for surfaces with $p_g>0$. Under the same assumptions, Theorem \ref{t4} of Section \ref{s1} provides an explicit expression for the cobordism series.

\subsection {A virtual Segre/Verlinde correspondence}
Let $X$ be a nonsingular projective surface, and let $\alpha \in K^0(X)$.
Using the usual (non-virtual) geometry of the
Hilbert scheme of points, we define
Segre and Verlinde series as follows:
\begin{eqnarray*} \mathsf S^{\mathsf{Hilb}}_{\alpha} &=&\sum_{n=0}^{\infty} q^n \int_{X^{[n]}} s(\alpha^{[n]})\, ,\\
\mathsf V^{\mathsf {Hilb}}_{\alpha} &=&\sum_{n=0}^{\infty} q^n\, \chi(X^{[n]}, \det \alpha^{[n]})
\,  .\end{eqnarray*}
Precise closed-form expressions for $\mathsf S^{\mathsf {Hilb}}_{\alpha}$
and $\mathsf V^{\mathsf {Hilb}}_{\alpha}$
are in general difficult to write down, but see \cite {EGL, Le, MOP, MOP3, V} for results and conjectures when $\alpha$ has small rank. Nevertheless, in the absence of explicit expressions,
a connection $$
\mathsf S^{\mathsf {Hilb}}_{{\alpha}}\, \, 
\leftrightarrow \,\,
\mathsf V^{\mathsf {Hilb}}_{\widetilde{\alpha}}
$$
after an explicit change of variables, 
for pairs $(\alpha, \widetilde \alpha)$ related in an explicit fashion,
was proposed in \cite {J, MOP2}.
The resulting {\it Segre/Verlinde correspondence}
is aligned with the larger conjectural framework of strange duality.  An extension to higher rank moduli spaces of sheaves was announced in \cite{GK}. 

In the virtual context, the Segre/Verlinde correspondence takes
a simpler form. We define the virtual (shifted)\footnote{The shift is needed to ensure a uniform statement in Theorem \ref{t2}.} Segre and Verlinde series \footnote{By contrast with the series considered in Section \ref{kthry}, the Segre/Verlinde series are given by algebraic functions; see Section \ref{s5} for details.} by
\begin{eqnarray*}
\mathsf S_{\alpha} (q)&=&q^{\beta. K_X}\sum_{n\in \mathbb Z} q^n \int_{\left[\mathsf{Quot}_{X}(\mathbb C^N, \beta, n)\right]^{\textrm{vir}}} s(\alpha^{[n]})\, ,\\
  \mathsf V_{\alpha} (q)&=&q^{\beta. K_X}\sum_{n\in \mathbb Z} q^n \,\chi^{\textrm{vir}}(\mathsf{Quot}_{X}(\mathbb C^N, \beta, n),\, \det \alpha^{[n]})\, .\end{eqnarray*} The more precise notation $\mathsf S_{N, \alpha}$ and $\mathsf V_{N, \alpha}$ also keeps track of the dimension of  $\mathbb{C}^N$.

\begin{theorem}\label{t2} The virtual Segre and Verlinde series match $$
\mathsf S_{N,\alpha}\left((-1)^{N}q\right) =
  \mathsf V_{N,\alpha}(q)$$ in the following three cases:
\begin{itemize}
\item [(i)] $X$ is a nonsingular projective surface, $\beta=0$, $N$ is arbitrary, 
\item [(ii)] $X$ is a nonsingular projective surface
  with $p_g>0$, $N=1$ and $\beta$ is arbitrary, 
\item [(iii)] $X$ is a relatively minimal elliptic surface, $\beta$ is supported on fibers, $N$ is arbitrary.
\end{itemize}
\end{theorem} 

However, the virtual Segre/Verlinde correspondence does not always hold. The simplest counterexample occurs for minimal surfaces of general type, $N=2$ and $\beta=K_X$, with discrepancy even at the first few terms. 

\subsection{Symmetry} For $\beta=0$, we show that the expressions for the Segre series $\mathsf S_{N, \,\alpha}$ are symmetric when the rank $r=\text{rk } \alpha$ and the dimension of
$\mathbb{C}^N$ are interchanged. 

\begin{theorem}\label{ss} Let $X$ be a nonsingular projective surface, and let $\alpha,
  \widetilde{\alpha} \in K^0(X)$ satisfy
  $$\mathsf{rk }\, \alpha=r,\quad \mathsf{rk }\, \widetilde \alpha=N\, ,\ \
  \quad\frac{\langle K_X, c_1(\alpha)\rangle }{\mathsf{rk }\, \alpha}=\frac{\langle K_X,c_1(\widetilde{\alpha})\rangle}{\mathsf{rk }\, {\widetilde \alpha}}\, .$$ Then, the following
  symmetry holds:
  $$\mathsf S_{N,\,{\alpha}}\left((-1)^{N}q\right)=\mathsf S_{r,\,\widetilde \alpha}\left((-1)^{r}q\right)\, .$$
\end{theorem}

\subsection{The virtual structure sheaf} 

Our calculations crucially rely on the simpler form of the virtual fundamental class of Quot schemes for curve classes of Seiberg-Witten length $N$. When $\beta=0$,
the following result was proven in \cite [Section 4.2.3]{OP}.  Let $C\in |K_X|$ be a nonsingular
canonical curve, and let $$\iota: \mathsf{Quot}_C(\mathbb C^N, n)\hookrightarrow \mathsf {Quot}_{X}(\mathbb C^N, n)$$ denote the natural inclusion. Then, the virtual class localizes on the canonical curve: 
\begin{equation}\label{localization}\left[\mathsf{Quot}_{X}(\mathbb C^N, n)\right]^{\vir}=(-1)^{n} \iota_{\star} \left[\mathsf {Quot}_C(\mathbb C^N, n)\right]\, .\end{equation} 

We provide a similar expression for the virtual structure sheaf.  
Let $\Theta=N_{C/X}$ denote the theta characteristic. Write $$\mathcal D_n=p^{\star}\det \left(\Theta^{[n]}\right)$$ for the pullback via the support morphism $$p: \mathsf {Quot}_{C}(\mathbb C^N, n)\to C^{[n]},\,\,\, Q\mapsto \text{supp }Q.$$ 
\begin{theorem} Assume $C\subset X$ is a nonsingular canonical curve. In rational $K$-theory, the virtual structure sheaf localizes on the canonical curve \label{tt}$$\mathcal O_{\mathsf{Quot}_X(\mathbb C^N, n)}^{\vir}=(-1)^{n} \,\iota_{\star}\, \mathcal D_n\, .$$ Furthermore, letting $\mathcal N^{\vir}$ denote the virtual normal bundle of the embedding $\iota$, the sheaf $\mathcal D_n$ is a square root of $\det \mathcal N^{\vir}$. 
\end{theorem} 
\subsection{The 8-fold equivalence}\label{sss}
Let $X$ be a nonsingular projective surface with a nonsingular canonical curve $C \subset X$.
Putting together Theorems \ref{t2} and \ref{ss}, we match the following $4$ virtual invariants of Quot schemes of $X$ in the $\beta=0$ case:
\begin{itemize}
\item [(i)] the Segre integrals $(-1)^{Nn} \int_{\left[ \mathsf{Quot}_X(\mathbb C^N, n)\right]^{\vir}} s(\alpha^{[n]})$
\item [(ii)] the Segre integrals $(-1)^{rn} \int_{\left[ \mathsf{Quot}_X(\mathbb C^r, n)\right]^{\vir}}  s(\widetilde{\alpha}^{[n]})$
\item [(iii)] the virtual Verlinde numbers $\chi^{\vir}( \mathsf{Quot}_X(\mathbb C^N, n), \,\,\,\det \alpha^{[n]})$
\item [(iv)] the virtual Verlinde numbers $\chi^{\vir}(\mathsf{Quot}_X(\mathbb C^r, n), \,\, \det \widetilde \alpha^{[n]}).$
\end{itemize} With the aid of Theorem \ref{tt} and equation \eqref{localization}, the same results can be stated over the canonical curve. For simplicity, let $\alpha=L^{\oplus r}, \widetilde {\alpha}=L^{\oplus N}$, where $L\to C$ is a line bundle. We match:
\begin{itemize}
\item [(i)$'$] the Segre integrals $(-1)^{Nn} \int_{ \mathsf{Quot}_C(\mathbb C^N, n)} s(L^{[n]})^r$
\item [(ii)$'$] the Segre integrals $(-1)^{rn} \int_{\mathsf{Quot}_C(\mathbb C^r, n)} s(L^{[n]})^{N}$
\item [(iii)$'$] the twisted Verlinde numbers $\chi\left( \mathsf{Quot}_C(\mathbb C^N, n), \,\,\,\left(\det L^{[n]}\right)^{r} \otimes \mathcal {D}_n\right)$
\item [(iv)$'$] the twisted Verlinde numbers $\chi\left(\mathsf{Quot}_C(\mathbb C^r, n), \,\, \left(\det L^{[n]}\right)^{N}\otimes \mathcal D_{n}\right).$
\end{itemize}

\noindent The diagram in Figure 1 summarizes the above 8-fold equivalence. 

\begin{figure}\vskip.2in
\begin{tikzcd}[back line/.style={densely dotted}, row sep=4.5em, column sep=-2.85em]
& \text{virtual Segre\,} \mathsf{Quot}_X(\mathbb C^N, n) \ar{dl}[swap, sloped,near start]{\text{cosection}} \ar{rr}[pos=3/7]{\text{Segre/Verlinde}} \ar[back line]{dd}[swap, sloped, above, pos=1/6]{\text{symmetry}} 
 & &\text{virtual Verlinde\,} \mathsf{Quot}_X(\mathbb C^N, n) \ar{dd}[sloped, above, pos=1/3]{\text{symmetry}} \ar{dl}[swap,sloped,near start]{\text{cosection}} \\
\text{Segre\,} \mathsf{Quot}_C(\mathbb C^N, n) \ar[crossing over]{rr}[pos=3/5]{\text{Segre/Verlinde}} \ar{dd}[sloped, above, pos=1/3]{{\text{symmetry}}} 
 & & \text{Verlinde }\mathsf{Quot}_C(\mathbb C^N, n) 
 \ar{dd}[swap, sloped, above, pos=1/6]{{\text{symmetry}}}  \\
& \text{virtual Segre }\mathsf{Quot}_X(\mathbb C^r, n) \ar[back line]{rr}[pos=2/5]{\text{Segre/Verlinde}} \ar[back line]{dl}[swap, sloped, near start] {\text{cosection}}
  & & \text{virtual Verlinde }\mathsf{Quot}_X(\mathbb C^r, n) \ar{dl} [swap, sloped, near start] {\text{cosection}}\\
\text{Segre\,} \mathsf{Quot}_C(\mathbb C^r, n) 
\ar{rr}[pos=3/5]{\text{Segre/Verlinde}} & & \text{Verlinde }\mathsf{Quot}_C(\mathbb C^r, n)\ar[crossing over, leftarrow]{uu}
\end{tikzcd}

\caption{The Segre/Verlinde series for curves/surfaces and their symmetry}

\vskip.2in
\end{figure}
\begin{question} Find a geometric interpretation of the above equalities. 
\end{question} 
It would be satisfying to see (i)$'$-(iv)$'$ as solutions to the same enumerative problem. The match between (iii)$'$ and (iv)$'$:
 $$\chi\left( \mathsf{Quot}_C(\mathbb C^N, n), \,\,\,\left(\det L^{[n]}\right)^{r} \otimes \mathcal {D}_n\right)=\chi\left(\mathsf{Quot}_C(\mathbb C^r, n), \,\, \left(\det L^{[n]}\right)^{N}\otimes \mathcal D_{n}\right)$$
 is reminiscent of strange duality for bundles over curves \cite{Bel, MO}. We will investigate these matters elsewhere. 

 \subsection{Plan of the paper}
  Section \ref{s1} is the most demanding computationally, but it plays a central role in our arguments. It generalizes the techniques of \cite {JOP, L, OP} and records structural results for the series of invariants we consider. It also contains a proof of Theorem \ref{p1} regarding the cobordism series. Section \ref{s3} shows the rationality of the $K$-theoretic series in several contexts: punctual quotients, curve classes of Seiberg-Witten length $N$, surfaces with $p_g=0$, and certain rank $1$ quotients. Section \ref{s4} establishes Theorem \ref{tt} regarding cosection localization for punctual quotients. Finally, Section \ref{s5} discusses the Segre/Verlinde correspondence and the symmetry of the Segre/Verlinde series, proving Theorems \ref{t2} and \ref{ss}.

\subsection{Acknowledgements} We thank A. Marian, M. Kool, T. Laarakker, A. Oblomkov, A. Okounkov, M. Savvas, S. Stark, R. Thomas and M. Zhang for
related discussions. We are especially grateful to A. Marian for conversations regarding the symmetries and equivalences in Section \ref{sss}. Our collaboration started at the {\it Algebraic Geometry and
  Moduli Zoominar} at ETH Z\"urich in the Spring of 2020.

N. A. was supported by the NSF through grant DMS-1902717.
D. J. was supported by SNF-200020-182181.
D. O. was supported by the NSF through grant DMS 1802228. R.P. was supported by
the Swiss National Science Foundation and
the European Research Council through
grants SNF-200020-182181, 
ERC-2017-AdG-786580-MACI. R.P. was also supported by SwissMAP.

The project has received funding from the European Research
Council (ERC) under the European Union Horizon 2020 Research and
Innovation Program (grant No. 786580).

\section{A general calculation and cobordism}\label{s1}

\subsection{Setup} While our interest here lies foremost in the series which appear in Theorems \ref{n2}, \ref{n1} and  \ref{p1} and \ref{t2}, it is useful to work more generally. 
Let $$f(x)\in F[\![x]\!]$$ be an invertible power series with coefficients in a field $F$. For our purposes, $F$ will be either the field of complex numbers or the field of rational functions in one or several variables. For a vector bundle $V$ over a complex projective scheme $S$, we consider  
the corresponding genus 
$$f(V)=\prod_{i}f(x_i)\in H^{\star}(S, F),$$ 
where $x_i$ are the Chern roots of $V$. This definition extends multiplicatively to $K$-theory.

Fix invertible power series $$f_1, \ldots, f_\ell, g\in F[\![x]\!],\quad g(0)=1,$$ and classes $$\alpha_1, \ldots, \alpha_\ell\in K^0(X) \text{ with }\text{rank }\alpha_s=r_s,\,\,\, 1\leq s\leq \ell.$$ As explained in Section \ref{kthry}, this yields $$\alpha_1^{[n]}, \ldots, \alpha_{\ell}^{[n]}\to \Quot_X(\C^N,\beta,n),$$ and we may consider the associated characteristic classes $$f_1(\alpha_1^{[n]}), \ldots, f_{\ell}(\alpha_{\ell}^{[n]})$$ in the cohomology of $\Quot_X(\C^N,\beta,n).$ Recall furthermore that the standard deformation-obstruction theory of the Quot scheme is governed by $$\text{Ext}^{0}(S, Q), \quad \text{Ext}^1(S, Q), \quad \text{Ext}^{\geq 2}(S, Q)=0.$$ Here, we keep the notation in \eqref {qqqq}, with $S, Q$ denoting the subsheaf and the quotient respectively. The last $\text{Ext}$ vanishes by Serre duality and the assumption that $\text{rank }Q=0$, see \cite {MOP}. In particular, the virtual tangent bundle is given by the alternating sum $$T^{\text{vir}}\Quot =\text{Ext}^{\bullet}(S, Q)=\text{Ext}^0(S, Q)-\text{Ext}^1(S, Q).$$

We consider the following general expression mixing the tautological classes and the virtual tangent bundle
\begin{equation}\label{eg}\mathsf {Z}_{X,N,\beta}^{f_1, \ldots, f_\ell ,g}(q\,|\,\alpha_1, \ldots, \alpha_\ell)=\sum_{n\in\mathbb Z} q^n \int_{[\Quot_X(\C^N,\beta,n)]^\vir} f_1(\alpha_1^{[n]})\cdots f_\ell(\alpha_\ell^{[n]}) \cdot g(T^\vir \Quot).\end{equation} As we will see below, the examples in Theorems \ref{n2}, \ref{n1}, \ref{p1}, \ref{t2} correspond to  
\begin{itemize}
\item [(i)] the virtual tautological series: $f_s(x)=1+y_s e^{x}$ and $g(x)=\frac{x}{1-e^{-x}},$  
\item [(ii)] the virtual Verlinde series: $f(x)=e^x$ and $g(x)=\frac{x}{1-e^{-x}},$ 
\item [(iii)] the virtual Segre series: $f(x)=\frac{1}{1+x}$ and $g(x)=1$,
\item [(iv)] the virtual cobordism class: $f(x)=1$ and $g$ arbitrary. 
\end{itemize}

An example covered by (iv), but not considered in further detail here is given by \begin{itemize}
\item [(v) ]the virtual elliptic genus: $f(x)=1$ and $g(x)=x\cdot \frac{\theta\left(\frac{x}{2\pi i}-z, \tau\right)}{\theta\left(\frac{x}{2\pi i}, \tau\right)}$, with $\theta$ denoting the Jacobi theta function.\footnote{The virtual elliptic genus has interesting automorphic properties for virtual Calabi-Yau manifolds. However, for $N=1$, $\beta=0$, the virtual canonical bundle $$K_{X^{[n]}}^{\text{vir}}=(\det TX^{[n]})^{\vee}\otimes \det \Obs= (K_X)_{(n)}\otimes \det \left((K_X)^{[n]}\right)^{\vee}=E^{\vee}$$ where $E$  is $-\frac{1}{2}$ of the exceptional divisor. To get weak Jacobi forms, the modified elliptic genus can be used: $$g(x)=x\cdot \frac{\theta\left(\frac{x}{2\pi i}-z, \tau\right)}{\theta\left(\frac{x}{2\pi i}, \tau\right)}\cdot \exp \left(-\frac{x}{2\pi i} \cdot \frac{\theta'}{\theta}(z, \tau)\right).$$ Alternatively, weak Jacobi forms can be obtained by twisting with pure classes coming from the symmetric product; these pair trivially with $K_{X^{[n]}}^{\text{vir}}$ and Theorem 5.4 of \cite {FG} applies.}
 \end{itemize}

\subsection {Calculation} \label{calcu}We compute \eqref{eg} explicitly for curve classes $\beta$ of Seiberg-Witten length $N$ following the approach in \cite{L, OP, JOP}. The strategy is to combine torus localization with universality arguments over the Hilbert scheme of points. Specifically, we let $\mathbb C^{\star}$ act with distinct weights $$w_1, \ldots, w_N$$ on the middle term of the exact sequence $$0\to S\to \mathbb C^N\otimes \mathcal O_X\to Q\to 0,$$ thus inducing an action on $\Quot_X(\C^N,\beta,n)$.
The series \eqref{eg} admits a natural equivariant lift. Writing $K$ for the canonical bundle of $X$ to ease notation, we show
\begin{theorem} For fixed $N,$ $r_1,\ldots, r_{\ell}$, there are universal series $$\mathsf A, \quad \mathsf B_s, \quad \mathsf U_{i}, \quad \mathsf V_{i, s}, \quad \mathsf W_{i,j}$$  such that for all nonsingular projective surfaces $X$, all curve classes $\beta$ of Seiberg-Witten length $N$, and all classes $\alpha_1, \ldots, \alpha_{\ell}$ of ranks $r_1, \ldots, r_{\ell}$, the equivariant series above takes the form
\begin{align}
\label{usus}
\mathsf {Z}_{X,N,\beta}^{f_1, \ldots, f_{\ell} ,g}&(q\,|\,\alpha_1, \ldots, \alpha_{\ell})=q^{-\beta.K}
\sum_{\substack{\beta=\beta_1+\cdots+\beta_N}}
\mathsf {SW}(\beta_1)\cdot\ldots\cdot \mathsf {SW}(\beta_N)\\& \cdot \left(\mathsf{A}^{K^2}\cdot \prod_{s=1}^{\ell}\mathsf{B}_s^{c_1(\alpha_s).K}\cdot \prod_{i=1}^N \mathsf{U}_i^{\,\beta_i.K} \cdot \prod_{\stackrel{1\leq i\leq N}{1\leq s\leq \ell}} \mathsf{V}_{i, s}^{\,\beta_i.c_1(\alpha_s)} \cdot \prod_{i<j} \mathsf{W}_{i,j}^{\, \beta_i.\beta_j}\right)\notag\,.
\end{align} The universal expressions appearing in the answer are explicitly recorded in equations \eqref{chvar} -- \eqref{seter} below. \end{theorem}

\noindent {\it Proof.} Under the above torus action, the $\mathbb C^{\star}$-fixed subsheaves split $$S=\bigoplus_{i=1}^{N} I_{Z_i}(-D_i),\quad [D_i]=\beta_i, \quad \text{length} (Z_i)=m_i$$ where $$\beta=\beta_1+\ldots+\beta_N, \quad m=m_1+\ldots+m_N,\quad m=n+\frac{1}{2} \sum_{i=1}^{N} \beta_i(\beta_i+K).$$ Thus, the torus fixed loci of $\Quot_X(\C^N,\beta,n)$ are isomorphic to products of Hilbert schemes of points and curves $$\mathsf F[\underline m, \underline \beta]=X^{[\underline m]}\times \Hilb_{\underline \beta}(X).$$ The underline notation 
$$\underline m=(m_1, \ldots, m_N), \quad \underline \beta=(\beta_1, \ldots, \beta_N)$$ is used to denote vectors of integers, curve classes etc. We set $$X^{[\underline m]}=X^{[m_1]}\times\ldots\times X^{[m_N]}, \quad \Hilb_{\underline \beta}(X)=\Hilb_{\beta_1}(X)\times \ldots \times  \Hilb_{\beta_N}(X).$$
 
Crucially, the arguments of Section 2 in \cite{L}, originally written for the integrals computing the virtual $\chi_{y}$-genus, apply verbatim to the more general expressions considered here. Indeed, via the virtual localization theorem in \cite {GP}, the integrals \eqref{usus} can be written as sums of contributions from the fixed loci $\mathsf F[\underline m, \underline \beta]$. Whenever $\beta$ is a curve class of Seiberg-Witten length $N$, nonvanishing contributions only arise from splittings satisfying \begin{equation}\label{vdzeroo}\beta_i\cdot (K-\beta_i)=0.\end{equation} In all other cases, the Seiberg-Witten invariants in the class $\beta_i$ vanish. As explained in \cite {L}, this in turn forces the vanishing of the corresponding localization terms. Furthermore, when \eqref{vdzeroo} is satisfied, the Hilbert schemes of curves $\mathsf {Hilb}_{\beta_i}(X)$ carry virtual fundamental classes of dimension zero and of length $\mathsf {SW}(\beta_i)$. Thus
\begin{equation}\label{spa}\left[{\mathsf F}[\underline m, \underline \beta]\right]^{\vir}=\mathsf {SW}(\beta_1)\cdots \mathsf {SW}(\beta_N)\cdot \left[X^{[\underline m]}\right]^{\vir}\times \left[\text{point}\right].\end{equation} The virtual class of $X^{[\underline{m}]}$ is determined by the obstruction bundle 
$$\text{Obs}=\sum_{i}\text{Obs}_i,\quad \text{Obs}_i=\left((K-\beta_i)^{[m_i]}\right)^{\vee}.$$
Note that $m=n+\beta\cdot K_X$. The equivariant series $\mathsf Z_{X, N, \beta}$ takes the following structure
$$\mathsf {Z}_{X,N,\beta}^{\underline f ,g}(q\,|\,\underline \alpha)=q^{-\beta.K_X}
\sum_{\substack{\beta=\beta_1+\cdots+\beta_N}}
\mathsf {SW}(\beta_1)\cdots \mathsf {SW}(\beta_N) \cdot \widehat{\mathsf {Z}}_{X,N, (\beta_1, \dots, \beta_N)}^{\underline f ,g} (q, \underline{w}\,|\, \underline \alpha),
$$ where we used the above convention to underline quantities understood as vectors.  The last term in the expression above arises via equivariant localization \cite{GP} and is given by
$$\widehat{\mathsf {Z}}_{X,N,\underline{\beta}}^{\underline f ,g} (q, \underline{w}\,|\, \underline \alpha)=\sum_{m_i\geq 0} q^{\sum m_i} \int_{X^{[\underline{m}]}} \prod_{s=1}^{\ell}f_s\Big(\iota^{\star}\alpha_s^{[n]}\Big) \cdot
\frac{\mathsf{e}}{g}\Big(\sum_{i=1}^{N} \Obs_i\Big) \cdot \frac{g}{\mathsf{e}}(N^\vir)
\cdot g\Big(\sum_{i=1}^{N} T_{X^{[m_i]}}\Big)$$
where $\mathsf e$ is the equivariant Euler class and $N^{\vir}$ denotes the virtual normal bundle. Furthermore,   
\begin{itemize}
\item[(i)] $\Obs_i =\left(\big(K-D_i)^{[m_i]}\right)^\vee$,\vskip.07in
\item [(ii)] $N^\vir =\sum_{i\neq j} N_{ij} [w_j-w_i], \quad \text{where} \quad N_{ij}=\textnormal{RHom}_\pi(I_{\mathcal Z_i}(-D_i), \mathcal O-I_{\mathcal Z_j}(-D_j))$,\vskip.07in
\item [(iii)] $\iota^{\star}\alpha^{[n]}=\sum_{i=1}^{N}\alpha_{\beta_i}^{[m_i]}[w_i], \,\,\,\,\,\, \text{with } \alpha_{\beta_i}^{[m_i]}=\textnormal{R}\pi_*(\mathcal Q_i\otimes q^*\alpha)=\alpha(-D_i)^{[m_i]}+\C^{\chi(\alpha |_{D_i})}.$\vskip.07in
\end{itemize}
 Here $\mathcal Z_i$ denotes the universal subscheme of $X^{[m_i]}\times X$, and the equivariant weights are recorded within brackets. 
Furthermore, $D_i\in \Hilb_{\beta_i}(X)$ are any representatives for the point classes appearing in \eqref{spa}. The integral does not depend on this choice.

The obstruction and normal bundles in (i) and (ii) are lifted from \cite {L, OP}. Item (iii) concerns the splitting of the tautological sheaves $\alpha^{[n]}$ to the product $X^{[\underline m]}$. This uses the splitting of the universal quotient $$\mathcal Q=\sum_{i=1}^N \mathcal Q_i[w_i],\,\,\,\,\,\,\, \mathcal Q_i=\mathcal O-\mathcal O(-D_i)\otimes \mathcal I_{\mathcal Z_i}.$$ We point out the following notational ambiguity in (iii): the tautogical sheaf $\alpha^{[n]}$ on the left is restricted from the Quot scheme to the fixed locus, but the tautological sheaves on the right are over the Hilbert scheme of $m_i$ points. The notation $\alpha_{\beta_i}^{[m_i]}$ records the twist by the curve class $\beta_i$ in the universal quotient. 

The expresssion $\widehat{Z}$ is a tautological integral over products of Hilbert schemes of points. Using the arguments of \cite {EGL}, it follows that $\widehat Z$ admits a multiplicative universal formula
\begin{align*}
\widehat{Z}_{X,N,(\beta_1, \dots, \beta_N)}^{\underline f,g}(q, \underline{w}\,|\,\underline \alpha)
&= \mathsf{A}^{K^2}\cdot \prod_{s=1}^{\ell}\mathsf{B}_s^{K.c_1(\alpha_s)}\cdot \prod_{1\leq s, t\leq \ell}\mathsf C_{s, t}^{c_1(\alpha_s). c_1(\alpha_t)}\cdot \prod_{s=1}^{\ell}\mathsf{D}_s^{c_2(\alpha_s)}\cdot \mathsf{E}^{\chi(\mathcal O_X)}\\
& \cdot \prod_{i=1}^N \mathsf{U}_i^{\,\beta_i.K} \cdot \prod_{\stackrel{1\leq i\leq N}{1\leq s\leq \ell}} \mathsf{V}_{i, s}^{\,\beta_i.c_1(\alpha_s)} \cdot \prod_{1\leq i<j\leq N} \mathsf{W}_{i,j}^{\, \beta_i.\beta_j}.
\end{align*}
The above universal series are independent of the choice of surface and curve classes $\beta_i$, though they do depend on $N, r_s, \underline {f}, g$. We omitted the Chern numbers $\beta_i^2$ because we will only use the formula for $\beta_i^2=\beta_i.K$.

\quad We next determine relations between the universal series by specializing to simpler geometries. When $\beta_1=\dots=\beta_N=K_X$, the obstruction bundle $$\Obs_i=\big(\mathcal O_X^{[m_i]}\big)^{\vee}$$ admits a trivial summand, hence its Euler class vanishes, for all $i$ with $m_i>0$. The case $m_i=0$ corresponds to the Hilbert scheme isomorphic to a point, and the only equivariant contributions come from (iii). 
Since for $D_i\in |K_X|$ we have $$\chi(\alpha_s|_{D_i})=-r_s \cdot K^2 + K\cdot c_1(\alpha_s),$$ the series $\widehat Z$ becomes 
$$\prod_{s=1}^{\ell}\prod_{i=1}^N \left(\frac{1}{f_s(w_i)^{r_s}}\right)^{K^2}\cdot f_s(w_i)^{K.c_1(\alpha_s)}.$$ This also equals $$\mathsf{A}^{K^2}\prod_{s=1}^{\ell}\mathsf{B}_s^{K.c_1(\alpha_s)}\prod_{1\leq s, t\leq \ell}\mathsf C_{s, t}^{c_1(\alpha_s). c_1(\alpha_t)}\prod_{s=1}^{\ell}\mathsf{D}_s^{c_2(\alpha_s)} \mathsf{E}^{\chi(\mathcal O_X)} \prod_{i=1}^N \mathsf{U}_i^{\,K^2} \prod_{\stackrel{1\leq i\leq N}{1\leq s\leq \ell}} \mathsf{V}_{i, s}^{\,K.c_1(\alpha_s)} \prod_{1\leq i<j\leq N} \mathsf{W}_{i,j}^{K^2}.$$ By the independence of the Chern numbers, we conclude that
$$\mathsf{C}_{s, t}\ =\ \mathsf{D}_s\ =\ \mathsf{E}\ =\ 1.$$
Therefore we have 
\begin{equation}
\label{us}\widehat{Z}_{X,N,(\beta_1, \dots, \beta_N)}^{\underline f,g}(q, \underline{w}\,|\,\underline \alpha)
= \mathsf{A}^{K^2}\cdot \prod_{s=1}^{\ell}\mathsf{B}_s^{K.c_1(\alpha_s)}\cdot \prod_{i=1}^N \mathsf{U}_i^{\,\beta_i.K} \cdot \prod_{\stackrel{1\leq i\leq N}{1\leq s\leq \ell}} \mathsf{V}_{i, s}^{\,\beta_i.c_1(\alpha_s)} \cdot \prod_{1\leq i<j\leq N} \mathsf{W}_{i,j}^{\, \beta_i.\beta_j}.
\end{equation}

\quad Next we specialize $X$ to a $K3$ surface blown up at a point, with exceptional divisor $E=K_X$. Let $I\sqcup J=[N]$ be a partition of the set $[N]=\{1, \dots, N\}$ into two disjoint parts. Define a vector of curve classes 
\[\underline{E}_{I\sqcup J}=(\beta_1, \dots, \beta_N) \quad \text{where} \quad 
    \begin{cases}
      \beta_i=E, & \text{if}\ i\in I \\
      \beta_j=0 , & \text{if}\ j\in J.
    \end{cases}
\]
By \eqref{us}, we have  
\[\widehat{Z}_{X,N,\underline{E}_{I\sqcup J}}^{\underline {f},g}(q, \underline{w}\,|\,\underline{\alpha})
= \mathsf{A}^{-1}\cdot \prod_{s=1}^{\ell}\mathsf{B}_s^{E.c_1(\alpha_s)}
\cdot \prod_{i\in I} \mathsf{U}_i^{\, -1} \cdot \prod_{\stackrel{i\in I}{1\leq s\leq \ell}} \mathsf{V}_{i, s}^{\, E.c_1(\alpha_s)} \cdot \prod_{\stackrel{i_1<i_2}{i_1, i_2\in I}} \mathsf{W}_{i_1,i_2}^{\, -1}.
\]
By varying the degree $d_s=E.c_1(\alpha_s)$ and the partition $I\sqcup J=[N]$, we determine all universal series. In fact, it will be enough to consider the cases when $|I|=0, 1$ or $ 2$. 

From now on, to ease the notation, we assume $\ell=1$ so that  
$$\widehat{\mathsf {Z}}_{X,N,\underline{\beta}}^{f ,g} (q, \underline{w}\,|\, \alpha)=\sum_{m_i\geq 0} q^{\sum m_i} \int_{X^{[\underline{m}]}} f\Big(\sum_{i=1}^{N} \alpha_{\beta_i}^{[m_i]} [w_i]\Big) \cdot
\frac{\mathsf{e}}{g}\Big(\sum_{i=1}^{N} \Obs_i\Big) \cdot \frac{g}{\mathsf{e}}(N^\vir)
\cdot g\Big(\sum_{i=1}^{N} T_{X^{[m_i]}}\Big).$$
We will make the convention to  temporarily denote indices in $I$, $J$ and $[N]$ by $i, j$, and $k$, respectively. Since
\begin{align*}\Obs_i&=\big((K_X-E)^{[m_i]}\big)^\vee=\big(\mathcal O_X^{[m_i]}\big)^\vee\,\, \text{if}\ i\in I\\
\Obs_j&=\big((K_X-0)^{[m_i]}\big)^\vee=\big(E^{[m_j]}\big)^\vee\,\,\,  \text{if}\ j\in J,
\end{align*}
all contributions with $m_i>0$, for some $i\in I$, must vanish. Therefore we assume that $m_i=0$ for all $i\in I$. For $j\in J$, the Euler class of $$\Obs_j=\big(E^{[m_j]}\big)^\vee$$ can be represented up to sign by the Hilbert scheme of $m_j$ points on $E=\mathbb P^1$:
\[\mathsf e(\Obs_j)=(-1)^{m_j} \left[E^{[m_j]}\right]=(-1)^{m_j} \left[\mathbb P^{m_j}\right].
\]
We obtain $\widehat{Z}_{X,N,\underline{E}_{I\sqcup J}}^{f,g}(q, \underline{w}\,|\,\alpha)$ equals 
$$\!\!\!\!\!\!\sum_{\forall j\in J, \ m_j\geq 0} \!\!\!\!\!\!(-q)^{\sum m_j} \int_{\prod\limits_j \mathbb P^{m_j}}\!
\iota^*\left(f\Big(\sum_{k=1}^{N} \alpha_{\beta_k}^{[m_k]} [w_k]\Big) \cdot \frac{g}{\mathsf e}(N^\vir) \cdot g\Big(\sum_j \left(T_{X^{[m_j]}}-\Obs_j\right)\Big)\right)
$$
where 
\[\iota:\textstyle \prod_{j} \mathbb P^{m_j}\hookrightarrow \textstyle\prod_{j} X^{[m_j]}.
\]
Let $d=E.c_1(\alpha).$ From the calculations of \cite[Section 3]{L} and \cite[Section 5]{OP} respectively, we have
\begin{align*}
\textnormal{(i)}&\,\, \iota^*\alpha_{\beta_i}^{[m_i=0]}[w_i]&&=\C^{\chi(\alpha |_{E})}[w_i]=\C^{r+d}[w_i],\\
\textnormal{(ii)}&\,\, \iota^* \alpha_{\beta_j=0}^{[m_j]}[w_j]&&=\mathcal O(-h_j)^{\oplus rm_j}[w_j]+\big(\mathcal O-\mathcal O(-h_j)\big)^{\oplus(r+d)}[w_j],\\
\textnormal{(iii)}&\,\, \iota^*(T_{X^{[m_j]}}-\Obs_j)&&= \mathcal O(-h_j)^{\oplus m_j}+\mathcal O(h_j)-\mathcal O \\ 
\textnormal{(iv)}&\,\, \iota^* N^\vir&&= \sum_{j_1\neq j_2} \Big(\mathcal O(-h_{j_2})^{\oplus m_{j_2}}+\mathcal O(h_{j_1})-\mathcal O(h_{j_1}-h_{j_2})\Big)[w_{j_2}-w_{j_1}]\\
& &&\ +\sum_{i, j} \Big(\mathcal O(-h_j)^{\oplus m_j}[w_j-w_i]+\mathcal O(h_j)[w_i-w_j]\Big).
\end{align*}
Here the $h_j$'s denote the pullbacks of the hyperplane classes from each of the factors. Integrating over the projective space means extracting the suitable coefficient of $h_j$; we denote this operation by square brackets. After carefully collecting all terms, we obtain that 
\begin{align*}\widehat{Z}_{X,N,\underline{E}_{I\sqcup J}}^{f,g}(q, \underline{w}\,|\,\alpha)
&=\sum_{m_j\geq 0} (-q)^{\sum m_j} \left[\textstyle\prod_{j\in J} h_j^{m_j}\right]\ \textstyle\prod_{j\in J} \Phi_j(h_j)^{m_j} \cdot \Psi_{I\sqcup J}(\{h_j\}_{j\in J})
\end{align*}
where 
\begin{align*}
\Phi_j(h_j)&=(-h_j)\cdot f(-h_j+w_j)^r\cdot \prod_{k=1}^N \frac{g(-h_j+w_j-w_k)}{-h_j+w_j-w_k}\\
\Psi_{I\sqcup J}&=\prod_{j\in J} h_j \cdot \prod_{k=1}^{N} f(w_k)^{r+d}\cdot \prod_{j\in J} \frac{1}{f(-h_j+w_j)^{r+d}}\cdot \prod_{\stackrel{ j\in J}{1\leq k\leq N}} \frac{g(h_j-w_j+w_k)}{h_j-w_j+w_k}\\&\cdot \prod_{j_1\neq j_2} \frac{h_{j_1}-w_{j_1}-h_{j_2}+w_{j_2}}{g(h_{j_1}-w_{j_1}-h_{j_2}+w_{j_2})}.
\end{align*}

We evaluate the above expression via the multivariable Lagrange-B\"urmann formula \cite{G}. Specifically, fix $\ell\geq 1$, and for each $1\leq j\leq \ell$, consider power series $\phi_j$ in the variable $x_j$, with nonzero constant terms. Set $$K={\prod_{j=1}^{\ell}\left(1-x_j\frac{d}{dx_j}\log \phi_j(x_j)\right)}.$$ Then, for all power series $\psi$ in $x_1, \ldots, x_{\ell},$ we have  
$$\sum_{m_j\geq 0} q^{m_1+\ldots+m_{\ell}}\left[x_1^{m_1}\cdots x_\ell^{m_{\ell}}\right]\ \phi_1(x_1)^{m_1} \cdot\ldots \cdot \phi_{\ell}(x_\ell)^{m_{\ell}} \cdot \psi(x_1, \ldots, x_{\ell})=\frac{\psi}{K}(x_1, \ldots, x_{\ell}).$$ On the right hand side, the $x_j$'s are the unique solutions to 
$$q=\frac{x_j}{\phi_j(x_j)},\quad x_j(q=0)=0.$$

Applying this to our situation, we obtain \[\widehat{Z}_{X,N,\underline{E}_{I\sqcup J}}^{f,g}(q, \underline{w}\,|\,\alpha)=
\frac{\Psi_{I\sqcup J}}{K_{I\sqcup J}}(\{h_j\}_{j\in J})\]
where $$K_{I\sqcup J}=\prod_{j\in J}\left(1-h_j\frac{d}{dh_j}\log \Phi_j(h_j)\right),$$
and for variables related by 
$$q=\frac{-h_j}{\Phi_j(h_j)},\quad h_j(q=0)=0.$$ 
In our case, we obtain
$$q=\frac{1}{f(-h_j+w_j)^r}\cdot \prod_{k=1}^{N} \frac{-h_j+w_j-w_k}{g(-h_j+w_j-w_k)} \quad \text{with}\quad h_j(q=0)=0$$
and 
$$K_{I\sqcup J}=\prod_{j\in J} h_j \cdot \prod_{j\in J} \left(r\cdot \frac{f'}{f}(-h_j+w_j)+\sum_{k=1}^{N} \Big(\frac{g'}{g}(-h_j+w_j-w_k)+\frac{1}{h_j-w_j+w_k}\Big)\right).$$
Using the additional change of variables 
$$H_j=h_j-w_j,$$
we can simplify this to 
$$\widehat{Z}_{X,N,\underline{E}_{I\sqcup J}}^{f,g}(q, \underline{w}\,|\,\alpha)=\frac{\prod\limits_{k=1}^{N} f(w_k)^{r+d}\cdot \prod\limits_{j\in J} \frac{1}{f(-H_j)^{r+d}}\cdot \prod\limits_{\stackrel{ j\in J}{1\leq k\leq N}}\frac{g(H_j+w_k)}{H_j+w_k}\cdot \prod\limits_{\stackrel{j_1\neq j_2}{j_1, j_2\in J}}\frac{H_{j_1}-H_{j_2}}{g(H_{j_1}-H_{j_2})} }{\prod\limits_{j\in J} \left(r\cdot \frac{f'}{f}(-H_j)+\sum\limits_{k=1}^{N} \Big(\frac{g'}{g}(-H_j-w_k)+\frac{1}{H_j+w_k}\Big)\right)},$$
where 
\begin{equation}\label{chvar}q=\frac{1}{f(-H_j)^r}\cdot \prod_{k=1}^{N} \frac{-H_j-w_k}{g(-H_j-w_k)} \quad \text{with}\quad H_j(q=0)=-w_j.\end{equation}

This determines all the universal series by comparing with 
\[\widehat{Z}_{X,N,\underline{E}_{I\sqcup J}}^{f,g}(q, \underline{w}\,|\,\alpha)
= \mathsf{A}^{-1}\cdot \mathsf{B}^{\, d}
\cdot \prod_{i\in I} \mathsf{U}_i^{\, -1} \cdot \prod_{i\in I} \mathsf{V}_i^{\, d} \cdot \prod_{\stackrel{i_1<i_2}{i_1, i_2\in I}} \mathsf{W}_{i_1,i_2}^{\, -1}.
\]
First, we find
\begin{equation}\label{bv}
\mathsf{B}=\prod_{k=1}^N \frac{f(w_k)}{f(-H_k)}, \quad \mathsf{V}_i=f(-H_i),
\end{equation}
by considering the terms with exponent $d$ and letting $I=\emptyset$ or $I=\{i\}$. Similarly, we see that 
\begin{eqnarray}
\label{A}
\mathsf{A}&=\prod\limits_{i=1}^{N}\frac{f(-H_i)^r}{f(w_i)^r}\cdot \prod\limits_{1\leq k,i\leq N}\,\frac{H_i+w_k}{g(H_i+w_k)}\,\cdot \prod\limits_{i_1\neq i_2}\frac{g(H_{i_1}-H_{i_2})} {H_{i_1}-H_{i_2}}
\cdot \\ &\cdot \prod\limits_{i=1}^{N} \left(r\cdot \frac{f'}{f}(-H_i)+\sum\limits_{k=1}^{N} \Big(\frac{g'}{g}(-H_i-w_k)+\frac{1}{H_i+w_k}\Big)\right)\notag
\end{eqnarray}
and
\begin{eqnarray}
\label{U}
\mathsf{U}_i&=\frac{1}{f(-H_i)^r}\cdot \prod\limits_{k=1}^{N}\frac{g(H_i+w_k)}{H_i+w_k}\cdot \prod\limits_{i'\neq i}\frac{H_{i'}-H_{i}}{g(H_{i'}-H_{i})}\cdot\frac{H_{i}-H_{i'}}{g(H_{i}-H_{i'})}\\ &\cdot \left(r\cdot \frac{f'}{f}(-H_i)+\sum\limits_{k=1}^{N} \Big(\frac{g'}{g}(-H_i-w_k)+\frac{1}{H_i+w_k}\Big)\right)^{-1}.\notag
\end{eqnarray}
by comparing the terms with exponent $(-1)$ and letting $I=\emptyset$ or $I=\{i\}$. In the above expressions, the indices $i, i'$ now range over $[N]$ because the partition $I\sqcup J=[N]$ is not involved. Finally, 
\begin{equation}\label{W}\mathsf{W}_{i_1, i_2}=\frac{g(H_{i_1}-H_{i_2})}{H_{i_1}-H_{i_2}}\cdot\frac{g(H_{i_2}-H_{i_1})}{H_{i_2}-H_{i_1}},
\end{equation}
by considering $I=\{i_1, i_2\}$ for $i_1<i_2$. 

When several functions $f_1, \ldots, f_\ell$ are involved, each instance of $f$ gets replaced by product contributions from each of the $f_s$'s. For instance, the change of variables becomes 
\begin{equation}\label{seter}q=\prod_{s=1}^{\ell}\frac{1}{f_s(-H_j)^{r_s}}\cdot \prod_{k=1}^{N} \frac{-H_j-w_k}{g(-H_j-w_k)}.\end{equation} This also affects the logarithmic derivatives $\frac{f'}{f}$ accordingly, turning them into sums of $\ell$ terms.\qed

\subsection{Cobordism and the virtual Pontryagin class} As an application of the calculations in the previous subsection, we consider the virtual cobordism series of the Hilbert scheme of curves and points: $$\mathsf Z^{\mathsf{cob}}_{X, \beta}=\sum_{n\in \mathbb Z} q^n \left[\mathsf{Quot}_{X}(\mathbb C^1, \beta, n)\right]^{\textnormal{vir}}_{\textnormal{cob}}\in \Omega^{\star}((q)).$$ 

For a proper scheme $S$ with a $2$-term perfect obstruction theory, virtual cobordism \cite {Sh} is encoded by all virtual Chern numbers $$\int_{[S]^{\vir}} c_{i_1}(T^{\vir}S)\cdots c_{i_k} (T^{\vir}S).$$ For bookkeeping, we set $$g(x)=\sum_{k=0}^{\infty} y_k x^k,\,\, y_0=1,$$ and note that the virtual cobordism class of $S$ is determined by $$[S]_{\textnormal{cob}}^{\vir}=\int_{[S]^{\vir}} g(T^{\vir}S)\in \mathbb C[y_0, y_1, \ldots].$$ This is consistent with the fact that $\Omega_{\star}$ is a polynomial ring in infinitely many variables. 

For surfaces with $p_g>0$, the series $$\mathsf Z^{\mathsf{cob}}_{X, \beta}=\sum_{n\in \mathbb Z} q^n \int_{\left[\mathsf{Quot}_{X}(\mathbb C^1, \beta, n)\right]^\vir} g(T^{\vir}\mathsf{Quot}_{X}(\mathbb C^1, \beta, n))$$ can be 
calculated by letting $N=1$ and $f=1$ in the formulas of Section \ref{calcu}, also setting $w=0$ in the answer. Then, $$\mathsf A=\mathsf U^{-1}=\frac{H}{g(H)}\cdot \left(\frac{g'}{g}(-H)+\frac{1}{H}\right),\,\,\,\mathsf B=\mathsf V=1$$ for the change of variables $$q=\frac{-H}{g(-H)}.$$ Thus $$\mathsf Z^{\mathsf{cob}}_{X, \beta}(q)=q^{-K\cdot \beta} \cdot \mathsf {SW}(\beta)\cdot \mathsf A^{K^2-K\cdot \beta}.$$ 

It is easy to see that $\mathsf Z^{\mathsf{cob}}_{X, \beta}$ is not given by a rational function. This fails for instance for the specialization $$y_1=0,\,\, y_2=y,\,\,\, y_k=0 \text{ for }k\geq 3\implies g(x)=1+yx^2$$ corresponding to the virtual Pontryagin class. In this case $$q=\frac{-H}{1+yH^2},\quad \mathsf A=\frac{1-yH^2}{(1+yH^2)^2}\implies \mathsf A=\frac{1}{2}\left(1+\sqrt{1-4q^2y}\right)\cdot \sqrt{1-4q^2y}.$$

We can however write down a different expression for the cobordism series. The cobordism ring $\Omega_{\star}$ is generated by the classes of projective spaces $\mathbb P^0, \mathbb P^1, \ldots.$ Consider the generating series $$\mathsf P=\sum_{n=0}^{\infty} q^n \left[\mathbb P^n\right]\in \Omega_{\star}((q)).$$ More generally, for each $\ell$, we consider the projective spaces $\mathbb P^n$ endowed with the nontrivial obstruction theory $\mathcal O(\ell).$  When $\ell\geq 1$ this corresponds to the class of a smooth degree $\ell$ hypersurface. We denote by $\left[\mathbb P^n\right]_{(\ell)}$ the associated virtual cobordism class and define $$\mathsf P_{\ell}=\sum_{n=0}^{\infty} q^n \left[\mathbb P^n\right]_{(\ell)}\in \Omega_{\star}((q)).$$
\begin{theorem} \label{t4}For surfaces with $p_g>0$, we have 
$$\mathsf Z^{\mathsf{cob}}_{X, \beta}=\mathsf {SW}(\beta)\cdot q^{-K^2} \cdot \left(-\frac{\mathsf P_{-1}}{\mathsf P^2} \right)^{K^2-K\cdot \beta}.$$\end{theorem}
\proof 
The series $\mathsf P$ corresponds under our conventions to $$\mathsf P=\sum_{n=0}^{\infty} q^n \int_{\mathbb P^n} g(T\mathbb P^n)=\sum_{n=0}^{\infty} q^n \int_{\mathbb P^n} g(h)^{n+1}=\sum_{n=0}^{\infty} q^n \cdot \left([h^n] \,g(h)^{n+1}\right)=\frac{dh}{dq}$$ by Lagrange-B\"urmann formula for $q=\frac{h}{g(h)}.$ Note that $h=-H$ in the notation used above.  Similarly,  \begin{eqnarray*}\mathsf P_{\ell}&=&\sum_{n=0}^{\infty} q^n \int_{\left[\mathbb P^n\right]_{(\ell)}} g(T^{\vir}\left[\mathbb P^n\right]_{(\ell)})=
\sum_{n=0}^{\infty} q^n \int_{\mathbb P^n} \mathsf e(\Obs)\cdot g(T\mathbb P^n-\Obs)\\&=&\sum_{n=0}^{\infty} q^n \int_{\mathbb P^n}\ell h\cdot \frac{g(h)^{n+1}}{g(\ell h)}=\sum_{n=0}^{\infty} q^n \cdot \left([h^n] \,\frac{\ell h\cdot g(h)}{g(\ell h)}\cdot g(h)^{n}\right)=\frac{\ell h}{g(\ell h)}\cdot \frac{dh}{dq},\end{eqnarray*} again using Lagrange-B\"urmann in the last equality. By direct calculation, we find $$\mathsf A=\frac{-h}{g(-h)}\cdot \left(\frac{g'}{g}(h)-\frac{1}{h}\right)=-\frac{1}{q}\cdot \frac{\mathsf P_{-1}}{\mathsf P^2},$$ completing the proof.\qed

\section{Rationality in K-theory and examples} \label{s3}

\subsection{Overview} We here establish Theorems \ref{n2} and \ref{n1} using the general calculations of Section \ref{calcu}. In addition, we justify Examples \ref{exa}, \ref{exc} and \ref{exd} in the Introduction. Finally, we consider surfaces with $p_g=0$, thereby proving Theorem \ref{n3}, and we discuss a calculation involving rank $1$ quotients. \vskip.1in

\subsection{Rationality} \label{srse}

{\it Proof of Theorem \ref{n1}.} The theorem concerns the shifted series $$\mathsf {\overline Z}=\sum_{n\in \mathbb Z} q^{n+\beta\cdot K} \chi^{\vir}(\mathsf {Quot}_X(\mathbb C^1, \beta, n), \wedge^{k_1}\alpha_1^{[n]}\otimes \ldots \otimes \wedge^{k_{\ell}} \alpha_\ell^{[n]}).$$ For formal variables $y_1, \ldots, y_{\ell}$, define $$\mathsf {\overline Z}^{\dagger}=\sum_{n\in \mathbb Z} q^{n+\beta\cdot K} \chi^{\vir}(\mathsf {Quot}_X(\mathbb C^1, \beta, n), \wedge_{y_1}\alpha_1^{[n]}\otimes \cdots \otimes \wedge_{y_{\ell}} \alpha_\ell^{[n]}).$$ Then $$\mathsf {\overline Z}=\frac{1}{k_1!\ldots k_{\ell}!} \cdot  \left (\frac{\partial}{\partial y_1}\right)^{k_{1}}\cdots \left (\frac{\partial}{\partial y_\ell}\right)^{k_{\ell}} \mathsf {\overline Z}^{\dagger}\bigg\vert_{\underline y=0}.$$ 

The series $\mathsf {\overline Z}^{\dagger}$ are special cases of the general expressions of Section \ref{s1}. This can be seen by invoking the Hirzebruch-Riemann-Roch theorem of \cite{CFK, FG}. For a scheme $S$ endowed with a $2$-term perfect obstruction theory, and any $K$-theory class $V\to S$, we have $$\chi^{\vir}(S, V)=\int_{\left[S\right]^{\vir}} \text{ch }(V)\cdot \text{Todd }(T^{\text{vir}}S).$$ In our context, we take $$S=\mathsf {Quot}_X(\mathbb C^1, \beta, n),\quad V=\wedge_{y_1}\alpha_1^{[n]}\otimes \ldots \otimes \wedge_{y_\ell} \alpha_\ell^{[n]}.$$ Thus,
in the notation of Section \ref{s1}, the series $\mathsf {\overline Z}^{\dagger}$ corresponds to functions $$f_s(x)=1+y_s e^{x},\quad 1\leq s\leq \ell,\quad g(x)=\frac{x}{1-e^{-x}}.$$  For $p_g>0$, we write via \eqref{usus}$$\mathsf {\overline Z}^{\dagger}=\mathsf {SW}(\beta) \cdot \mathsf A^{K^2-\beta\cdot K}\cdot \prod_{s=1}^{\ell} \mathsf B_s^{\,c_1(\alpha_s).K} \cdot \prod_{s=1}^{\ell} \mathsf V_s^{\,c_1(\alpha_s)\cdot \beta},$$ where we used $\mathsf A=\mathsf U^{-1}.$   
The change of variables \eqref{chvar} becomes $$q=(1-e^{H})\prod_{s=1}^{\ell} (1+y_s e^{-H})^{-r_s}.$$ For convenience, set $e^{H}=1-z$ so that $$q=z\cdot \prod_{s=1}^{\ell} \left(\frac{1-z}{1-z+y_s}\right)^{r_s}.$$ We regard $z$ as a function of $q, y_1, \ldots, y_\ell$. From formulas \eqref{bv}--\eqref{U} of Section \ref{calcu} we find \begin{align*}\mathsf A&=\mathsf U^{-1}=\prod_{s=1}^{\ell} \left(1+\frac{y_s}{1-z}\right)^{r_s}\cdot \prod_{s=1}^{\ell}(1+y_s)^{-r_s} \cdot \left(\sum_{s=1}^{\ell} r_s \frac{y_s}{1+y_s-z}\cdot \frac{-z}{1-z}+1\right)\\\mathsf B_s&=\frac{(1+y_s)(1-z)}{1+y_s-z}\\\mathsf V_s&=1+\frac{y_s}{1-z}.
\end{align*} We immediately see that $$\mathsf A\bigg\vert_{\underline{y}=0}=\mathsf B_s\bigg\vert_{\underline{y}=0}=\mathsf V_s\bigg\vert_{\underline{y}=0}=1.$$ 

Let us first assume $\ell=1$ for simplicity, also writing $y_1=y$ and $r_1=r$.  Note that $z\bigg\vert_{y=0}=q.$ The function $z$ is differentiable in $y$. We claim that for $k\geq 1$, the derivatives take the form $$\frac{\partial^k z}{\partial y^k}\bigg\vert_{y=0}=\frac{P_k(q)}{(1-q)^{2k-1}}$$ for certain polynomials $P_k$.\footnote{A more precise calculation similar to that in Lemma \ref{llll} below shows that the derivative equals $$(k-1)!\sum_{n=1}^{\infty} r\binom{rn-1}{k-1}\binom{n+k-2}{k-1} q^n.$$ The coefficient of $q^n$ is a polynomial in $n$ of degree $2k-2$, from where the claim follows.} This can be seen by induction on $k$ starting from the equation $$F(y, z)=z\cdot \left(\frac{1-z}{1-z+y}\right)^{r}=q.$$ For instance, differentiating once we obtain $$\frac{\partial{F}}{{\partial y}}+\frac{\partial{F}}{{\partial z}}\cdot \frac{\partial z}{\partial y}=0\implies \frac{\partial z}{\partial y}\bigg\vert_{y=0}=\frac{rq}{1-q}.$$ The case of arbitrary $k$ is notationally more involved, but it is useful to spell out the details in order to bound the order of the poles. Key to the argument is the general expression for the $k$th derivative of $F(y, z)$ with respect to $y$. It takes the form $$\sum \mathsf C\left(p \,| \,m_1, \ldots, m_t\right)\cdot \left(\frac{\partial}{\partial z}\right)^{t}\,\left(\frac{\partial}{\partial y}\right)^{p}\, F \cdot \frac{\partial^{m_1} z}{\partial y^{m_1}}\cdots \frac{\partial^{m_t} z}{\partial y^{m_t}}=0.$$ In the above expression, the length of the tuple $(m_1, \ldots, m_t)$ must equal the number of $z$ derivatives, namely $t$, and of course the total number of $y$ derivatives equals $$p+m_1+\ldots+m_t=k.$$ The derivative $\frac{\partial^k z}{\partial y^k}$ appears in the form $$\frac{\partial F}{\partial z} \cdot \frac{\partial^k z}{\partial y^k} + \text{ lower order terms}=0.$$ (The leading coefficient $\mathsf C(0\,|\, k)$ can be seen to be $1$ by induction.) The lower order terms involve lower derivatives of $z$ with respect to $y$ and derivatives of $F$. By direct calculation we find $\frac{\partial F}{\partial z}\bigg\vert_{y=0}=1$ so we can solve for the derivative $\frac{\partial^k z}{\partial y^k}\bigg\vert_{y=0}$ in terms of the lower order terms. In fact, using the explicit expression for $F$, all derivatives take the form $$ \left(\frac{\partial}{\partial z}\right)^{t}\,\left(\frac{\partial}{\partial y}\right)^{p}\, F\bigg\vert_{y=0} =\frac{F_{p, t}(q)}{(1-q)^{p+t}}$$ for certain polynomials $F_{p, t}$. Inductively, all terms $\frac{\partial^{m_1} z}{\partial y^{m_1}}\cdots \frac{\partial^{m_t} z}{\partial y^{m_t}}\bigg\vert_{y=0}$ are rational functions in $q$ with pole at $q=1$ of order less or equal to $(2m_1-1)+\ldots+(2m_t-1)$. Solving for $\frac{\partial^k z}{\partial y^k}\bigg\vert_{y=0}$ we obtain a rational function with pole of order at most $$(p+t)+(2m_1-1)+\ldots+(2m_t-1)=p+2m_1+\ldots+2m_t\leq 2k-1.$$ The only exception is $p=0$ which requires $t\geq 2$ (since $t=1$ gives the leading term), but then by direct calculation $$ \left(\frac{\partial}{\partial z}\right)^{t}F\bigg\lvert_{y=0}=0.$$ This completes the proof of the claim.
 
Now let $G(y, z)$ be one of the factors that appears in the product expressions of $\mathsf A,$ $\mathsf B,$ $\mathsf V$ or their inverses $\mathsf A^{-1},$ $\mathsf B^{-1},$ $\mathsf V^{-1}$. These terms are fractions in $y$ and $z$. We claim that  \begin{equation}\label{pro}\frac{\partial^k G}{\partial y^k}\bigg\vert_{y=0}\text{ is a rational function in } q \text { with pole of order at most }2k\text{ at }q=1.\end{equation} Indeed, the derivative is given by the same formula
$$\sum \mathsf C\left(p \,|\, m_1, \ldots, m_t\right)\cdot \left(\frac{\partial}{\partial z}\right)^{t}\,\left(\frac{\partial}{\partial y}\right)^{p}\, G \cdot \frac{\partial^{m_1} z}{\partial y^{m_1}}\cdots \frac{\partial^{m_t} z}{\partial y^{m_t}}.$$ In all cases at hand, we can compute 
$$\left(\frac{\partial}{\partial z}\right)^{t}\,\left(\frac{\partial}{\partial y}\right)^{p}\, G\bigg\vert_{y=0}$$ directly and observe that the answers are rational functions in $q$ with poles of order at most $2p+t$ at $q=1$. Thus, at $y=0$, the derivative in \eqref{pro} is rational with poles of order at most $$(2p+t)+(2m_1-1)+\ldots+(2m_t-1)= 2k.$$

Property \eqref{pro} is preserved under taking products. We therefore conclude that $\mathsf {\overline Z}^{\dagger}$ satisfies \eqref{pro} as well. Equivalently, $\mathsf {\overline Z}$ is rational in $q$ with pole of order at most $2k$ at $q=1$, as claimed. 

Finally, the case $\ell>1$ is similar. The above statements hold true, their proof requiring only more diligent bookkeeping via multiindices, but no new ideas. \qed
\vskip.1in

{\it Proof of Theorem \ref{n2}.}  For simplicity, first take $\ell=1$. We seek to prove the rationality of the series $$\mathsf {Z}_{X, N, \beta}(\alpha\, |\, k)=\sum_{n\in \mathbb Z} q^n \chi^{\textrm{vir}}\left(\mathsf{Quot}_{X}(\mathbb C^N, \beta, n),\wedge^k \alpha^{[n]}\right).$$ The proof exploits the symmetric structure of the expression \eqref{usus}, and also relies on the rewriting of the change of variables \eqref{chvar} and of the functions \eqref{bv} -- \eqref{W} in terms of rational functions, in the case at hand. 

Specifically, as in the proof of Theorem \ref{n1} above, set $$\mathsf Z^{\dagger}=\sum_{n\in \mathbb Z} q^n \chi^{\textrm{vir}}\left(\mathsf{Quot}_{X}(\mathbb C^N, \beta, n),\wedge_y \alpha^{[n]}\right)\implies \mathsf {Z}=\frac{1}{k!}\frac{\partial^k \mathsf Z^{\dagger}}{\partial y^k} \bigg\rvert_{y=0}.$$ For Seiberg-Witten classes $\beta$ of length $N$, by \eqref{usus} we have $$\mathsf Z^{\dagger}=q^{-\beta\cdot K_X} \sum_{\beta_1+\ldots+\beta_N=\beta} \mathsf {SW}(\beta_1)\cdots \mathsf{SW}(\beta_N) \cdot \mathsf{A}^{K^2}\cdot \mathsf{B}^{K.c_1(\alpha)}
\cdot \prod_{i=1}^N \mathsf{U}_i^{\,\beta_i.K} \cdot \prod_{i=1}^N \mathsf{V}_i^{\,\beta_i.c_1(\alpha)} \cdot \prod_{i<j} \mathsf{W}_{i,j}^{\, \beta_i.\beta_j}.$$ 
For every decomposition $\beta=\beta_1+\ldots+\beta_N$, we consider the symmetrized expression $$\mathsf C=\frac{1}{|\text{Stab}|}\sum_{\sigma\in S_N} \prod_{i=1}^N \mathsf{U}_{i}^{\,\beta_{\sigma(i)}.K} \cdot \prod_{i=1}^N \mathsf{V}_{i}^{\,\beta_{\sigma(i)}.c_1(\alpha)} \cdot \prod_{i<j} \mathsf{W}_{i, j}^{\, \beta_{\sigma(i)}.\beta_{\sigma(j)}}$$ where the sum is over the orbit of the symmetric group $S_N$ permuting the classes $\beta_i$. We will drop the dependence of $\mathsf C$ on the classes $\beta_i$ for simplicity. 

We use the results of Section \ref{calcu} for the functions $$f(x)=1+ye^x,\quad g(x)=\frac{x}{1-e^{-x}}.$$ To simplify the formulas, we further set $e^{-H_j}=1-z_j,\,\, e^{w_k}=1+t_k.$ The change of variables \eqref{chvar} becomes \begin{equation}\label{changev}q=\left(1+{y}(1-z_j)\right)^{-r} \cdot \prod_{k=1}^{N}\left(\frac{-z_j-t_k}{1-z_j}\right),\quad z_j(0)=-t_j.\end{equation} We record \begin{align*}
\mathsf A&=\prod_{i=1}^{N}\left(\frac{1+y(1-z_i)}{1+y(1+t_i)}\right)^{r} \prod_{1\leq i, k\leq N} \frac{z_i+t_k}{1+t_k}\prod_{i\neq j} \frac{1-z_j}{z_i-z_j} \prod_{i=1}^{N} \left(\frac{ry(1-z_i)}{1+y(1-z_i)} + \sum_{k=1}^{N} \frac{1+t_k}{z_i+t_k}\right)\\
\mathsf B&=\prod_{i=1}^{N}\frac{1+y(1+t_i)}{1+y(1-z_i)}\\
\mathsf U_i&=\frac{1}{(1+y(1-z_i))^r} \prod_{k=1}^{N}\frac{1+t_k}{z_i+t_k} \prod_{i'\neq i} \left(\frac{z_i-z_{i'}}{1-z_{i'}} \frac{z_{i'}-z_i}{1-z_i}\right) \left( \frac{ry(1-z_i)}{1+y(1-z_i)} + \sum_{k=1}^{N} \frac{1+t_k}{z_i+t_k}\right)^{-1}   \\
\mathsf V_i&=1+y(1-z_i)\\
\mathsf W_{i_1, i_2}&=\frac{1-z_{i_1}}{z_{i_2}-z_{i_1}}\frac{1-z_{i_2}}{z_{i_1}-z_{i_2}}.
\end{align*}

To show that $$\mathsf Z=\frac{1}{k!} \frac{\partial^k \mathsf Z^{\dagger}}{\partial y^k} \bigg\rvert_{y=0,\,\underline{t}=0}$$ is a rational function in $q$, it suffices to show that the derivatives  
 \begin{equation}\label{ddd}\frac{\partial^k \mathsf A}{\partial y^k}\bigg\rvert_{y=0, \,\underline{t}=0},\quad \frac{\partial^k \mathsf B}{\partial y^k}\bigg\rvert_{y=0, \,\underline{t}=0},\quad \quad \frac{\partial^k \mathsf C}{\partial y^k}\bigg\rvert_{y=0, \,\underline{t}=0},\end{equation}
have the same property. We do not have closed form expressions for these derivatives, but a qualitative argument suffices. 

Property \eqref{ddd} is in fact valid for any rational function $\mathsf R$ in the $\{z_i\}$, $\{t_k\}$ and $y$, which is symmetric in the $z$'s. 
The functions $\mathsf A, \mathsf B, \mathsf C$ are of this type. This is clear for $\mathsf A$ and $\mathsf B$ from the explicit expressions, and it also holds for $\mathsf C$ because we sum over the orbit of the $S_N$-action. 

The assertion above was proven in \cite {JOP}, Section $2.2$, and we only indicate the main points. Indeed, ${\mathsf R}\bigg\vert_{y=0,\,\, \underline{t}=0}$ is a rational function in $q$ since it can be expressed in terms of the elementary symmetric functions in the $z_i$'s. These symmetric functions are continuous in $y$ and $\underline{t}$. At $y=0,\, \underline{t}=0$, they are rational functions in $q$ since the change of variables becomes $$z^{N}=q(z-1)^N.$$ The derivative $$\frac{\partial \mathsf R}{\partial y} = \sum_{i=1}^{N} \frac{\partial \mathsf R}{\partial z_i} \cdot \frac{\partial z_i}{\partial y} + \partial_{y} \mathsf R$$ is 
also given by a rational function in  $\{z_i\}$, $\{t_k\}$ and $y$, symmetric in the $\{z_i\}$'s. 
The second term $\partial_{y} \mathsf R$ is calculated keeping $z$ fixed, and is symmetric in the $z$'s since $\mathsf R$ is. By implicit differentiation $$F(y, \,z, \,\underline{t})=(1+y(1-z))^{-r}\cdot \prod_{k=1}^{N} \frac{-z-t_k}{1-z}=q\implies \frac{\partial z}{\partial y} = \mathsf S(y, \,z, \,\underline{t})$$ for an explicit rational function $\mathsf S$.  Using that $\mathsf R$ is symmetric and the formula for the implicit derivatives, it follows that transposing $z_i$ and $z_j$ turns $$\frac{\partial \mathsf R}{\partial z_i} \text{ into  }\frac{\partial \mathsf R}{\partial z_j} \quad \text{and} \quad \frac{\partial z_i}{\partial y} \text{ into } \frac{\partial z_j}{\partial y}.$$ Since we sum over all $i$'s, $\frac{\partial \mathsf R}{\partial y}$ is symmetric in the $z_i$'s. 

Inductively, it follows that all higher derivatives of $ \frac{\partial^k \mathsf R}{\partial y^k}$ are rational functions of  $\{z_i\}$, $\{t_k\}$ and $y$, symmetric in the $\{z_i\}$. Their values $ \frac{\partial^k \mathsf R}{\partial y^k}\bigg\vert_{y=0,\,\underline{t}=0}$ are therefore rational in $q$. 

Finally, introducing additional formal variables $y_1, \ldots, y_\ell$ and running the argument above, we obtain the rationality of the series $$\sum_{n\in \mathbb Z} q^n\chi^{\vir}\left(\mathsf {Quot}_{X} (\mathbb C^N, \beta, n), \wedge^{k_1} \alpha_1^{[n]}\otimes \ldots \otimes \wedge^{k_{\ell}} \alpha_{\ell}^{[n]}\right).\quad{\qed}$$
\begin{remark} The same arguments prove the rationality of the generating series of descendent invariants defined in \cite {JOP}, for all curve classes $\beta$ of Seiberg-Witten length $N$. This partially confirms Conjecture $2$ in the work cited. 
\end{remark}

\subsection{Examples}\label{workex} Equations \eqref{taut1} and \eqref{taut2} in the Introduction are obtained by executing the above proof in more precise detail. We continue to use the notation $$\mathsf {Z}=\sum_{n\in \mathbb Z} q^n \chi^{\textrm{vir}}\left(\mathsf{Quot}_{X}(\mathbb C^N, n), \alpha^{[n]}\right),\,\,\, \mathsf {Z}^{\dagger}=\sum_{n\in \mathbb Z} q^n \chi^{\textrm{vir}}\left(\mathsf{Quot}_{X}(\mathbb C^N, n), \wedge_{ y } \alpha^{[n]}\right).$$ We have $$\mathsf {Z}=\frac{\partial \mathsf Z^{\dagger}}{\partial y} \bigg\rvert_{y=0},\,\,\,\mathsf Z^{\dagger}=\mathsf A^{K^2}\cdot \mathsf B^{c_1(\alpha).K}.$$ The functions $\mathsf A, \mathsf B$ were recorded in the proof of Theorem \ref{n2}. We claim that \begin{itemize}
\item [(i)] $\mathsf A\bigg\rvert_{y=0, \,\underline{t}=0}=\mathsf B\bigg\rvert_{y=0, \,\underline{t}=0}=1$ \vskip.05in
\item [(ii)] $\frac{\partial \mathsf A}{\partial y}\bigg\rvert_{y=0, \,\underline{t}=0}=-Nr\cdot \frac{q^2}{(1-q)^2},\quad \frac{\partial \mathsf B}{\partial y}\bigg\rvert_{y=0, \,\underline{t}=0}=-N\cdot \frac{q}{1-q}.$
\vskip.05in
\end{itemize}
From here, \eqref{taut1} and \eqref{taut2} follow immediately since
$$\mathsf Z=\frac{\partial \mathsf Z^{\dagger}}{\partial y} \bigg\rvert_{y=0}=\frac{\partial}{\partial y} \left(\mathsf A^{K^2}\cdot \mathsf B^{c_1(\alpha).K}\right)\bigg\rvert_{y=0}=- N r \cdot K^2\cdot \frac{q^2}{(1-q)^2} - N\cdot (c_1(\alpha)\cdot K)\cdot \frac{q}{1-q}.
$$

Item (i) is clear due to equation \eqref{taut0}. This can also be seen directly by setting $$y=0, \,\underline{t}=0$$ in the explicit formulas. The fact that these substitutions make sense is a consequence of the arguments in Section $2.2$ of \cite {JOP}: the elementary symmetric functions in the roots $z_1, \ldots, z_N$ are continuous in the parameters. 
When $y=0,\,\underline{t}=0$, the equation giving the change of variables becomes $$z^N-q(z-1)^N=0$$ with roots $\mathsf{z}_i$. It is then clear that $\mathsf B\bigg\vert_{y=0,\, \underline{t}=0}=1$ and $$\mathsf A{\bigg\rvert_{y=0,\, \underline{t}=0}}=\prod_{i=1}^{N} {\mathsf{z}_i^N}\cdot \prod_{i\neq j} \frac{1-\mathsf{z}_j}{\mathsf{z}_i-\mathsf{z}_j}\cdot \prod_{i=1}^{N} \frac{N}{\mathsf{z}_i}=1.$$

Item (ii) is similar. It is more convenient to work with the logarithmic derivative. Set $\mathsf B^{\dagger}=\log \mathsf B$. Using (i), we have \begin{align*}\frac{\partial \mathsf B}{\partial y}\bigg\rvert_{y=0, \,\underline{t}=0}=\frac{\partial \mathsf B^{\dagger}}{\partial y}\bigg\rvert_{y=0, \,\underline{t}=0} &=\sum_{i=1}^{N}\frac{\partial}{\partial y} \left(\log (1+y)-\log (1+y(1-z_i))\right)\bigg\rvert_{y=0} \\&=\sum_{i=1}^{N}\mathsf{z}_i=-N\cdot \frac{q}{1-q}.\end{align*} We record here that \begin{equation}\label{vie}\sum_{i=1}^{N} \mathsf {z}_i = -N\frac{q}{1-q},\quad\,\,\sum_{i<j} \mathsf{z}_i \mathsf{z}_j =-\binom{N}{2}\frac{q}{1-q}.\end{equation} The calculation involving $\mathsf A$ is slightly more involved and requires the implicit derivatives
$$\frac{\partial z_i}{\partial y}\bigg\rvert_{y=0, \,\underline{t}=0}=\frac{r}{N} \mathsf z_i (1-\mathsf z_i)^2.$$ These are found by differentiating \eqref{changev}. Writing $\mathsf A^\dagger=\log \mathsf A$, we obtain \begin{align*}\frac{\partial  \mathsf A}{\partial y}\bigg\rvert_{y=0, \,\underline{t}=0}&=\frac{\partial  \mathsf A^{\dagger}}{\partial y}\bigg\rvert_{y=0, \,\underline{t}=0}=\sum_{i=1}^{N} \frac{\partial  \mathsf A^{\dagger}}{\partial z_i}\bigg\rvert_{y=0, \,\underline{t}=0}\cdot \frac{\partial z_i}{\partial y}\bigg\rvert_{y=0, \,\underline{t}=0}+\frac{\partial  \mathsf A^{\dagger}}{\partial y}\bigg\rvert_{y=0, \,\underline{t}=0}.\end{align*} The last $y$-derivative is understood as keeping the $z$'s fixed. The answer is symmetric in the $z$'s and setting $y=\underline{t}=0$ is therefore allowed. Direct calculation gives 
$$\frac{\partial  \mathsf A^{\dagger}}{\partial y}\bigg\rvert_{y=0, \,\underline{t}=0}=-r\sum_{i=1}^{N} \mathsf{z}_i + \frac{r}{N} \sum_{i=1}^{N}\mathsf{z}_i(1-\mathsf{z}_i)=-Nr\cdot \frac{q^2}{(1-q)^2}.$$ Similarly, the following sum is also symmetric in the roots, and by direct computation we find \begin{align*}\sum_{i=1}^{N} \frac{\partial  \mathsf A^{\dagger}}{\partial z_i}\bigg\rvert_{y=0, \,\underline{t}=0}\cdot \frac{\partial z_i}{\partial y}\bigg\rvert_{y=0, \,\underline{t}=0}&=\frac{r}{N}\cdot \sum_{i=1}^{N} \frac{\partial  \mathsf A^{\dagger}}{\partial z_i}\bigg\rvert_{y=0, \,\underline{t}=0}\cdot \mathsf{z}_i(1-\mathsf{z}_i)^2\\&=\frac{r}{N}\left((N-1)\sum_{i} \mathsf{z}_i - 2 \sum_{i<j} \mathsf{z}_i \mathsf{z}_j\right)=0,\end{align*} using \eqref{vie}. This completes the argument.\vskip.1in

Finally, equations \eqref{taut3} and \eqref{taut4} of the Introduction require not only the series $\mathsf A$, $\mathsf B$ studied in (i)--(ii), but also $\mathsf U_i,$ $\mathsf V_i$ and their first derivatives. For $\text{rank }\alpha=0$, we find 
\begin{itemize}
\item [(iii)] $\mathsf U_i\bigg\rvert_{y=0, \,\underline{t}=0}= -\frac{Nq}{1-q}\cdot \frac{1}{\mathsf z_i},\quad \mathsf V_i\bigg\rvert_{y=0, \,\underline{t}=0}=1$ \vskip.05in
\item [(iv)] $\frac{\partial \mathsf U_i}{\partial y}\bigg\rvert_{y=0, \,\underline{t}=0}=0,\quad \frac{\partial \mathsf V_i}{\partial y}\bigg\rvert_{y=0, \,\underline{t}=0}=1-\mathsf z_i.$
\end{itemize} 
The expressions in Examples \ref{exc} and \ref{exd} in the Introduction follow by substituting (i) -- (iv) in \eqref{usus}. Items (i) -- (iv) and equation \eqref{W} suffice to treat all classes $\beta=\ell K_X$ for $0\leq \ell \leq N$, but the formulas are more cumbersome. 

\subsection{Surfaces with $p_g=0$} We next give a proof of Theorem \ref{n3}. Let $X$ be a simply connected surface with $p_g=0$ and set $N=1$. We show that the shifted series $$\overline{\mathsf Z}=q^{\frac{1}{2}\beta(\beta+K_X)}\cdot \sum_{n\in \mathbb Z} q^{n} \chi^{\vir}\left(\mathsf {Quot}_{X} (\mathbb C^1, \beta, n), \,\alpha^{[n]}\right)$$ is given by a rational function with poles only at $q=1$. This case reduces to the calculation of Riemann-Roch numbers of tautological sheaves over the (non-virtual) Hilbert scheme of points. 

We have the following identification $$\mathsf {Quot}_{X}(\mathbb C^1, \beta, n)\simeq X^{[m]}\times \mathbb P$$ where $\mathbb P=|\beta|$ and $m=n+\frac{\beta(\beta+ K_X)}{2}$. When $p_g=0$, the obstruction bundle was found in \cite {OP}: $$\text{Obs}= H^1(\beta)\otimes \zeta + \left(M^{[m]}\right)^{\vee}\otimes \zeta,$$ where $\zeta=\mathcal O_{\mathbb P}(1)$ denotes the tautological bundle, and $M=K_X-\beta$. Thus $$\mathcal O^{\vir}_{\mathsf {Quot}}= \wedge_{-1} \text{Obs}^{\vee}=\left(\wedge_{-1} \zeta^{-1} \right)^{h^1(\beta)}\otimes \wedge_{-1} \left(M^{[m]}\otimes \zeta^{-1}\right).$$ The term 
$$\wedge_{-1} \zeta^{-1}=\mathcal O-\zeta^{-1}$$ is supported on hyperplanes, and its powers can be used to cut down the dimension of the projective space. Hence $$\mathcal O^{\vir}_{\mathsf {Quot}}=\iota_{\star} \wedge_{-1} \left(M^{[m]}\otimes \zeta^{-1}\right)$$ is supported on $X^{[m]}\times \mathbb P(V)$ where $$\dim V = h^0(\beta)-h^1(\beta)=\chi(\beta)=\nu+1,\,\, \text{ assuming }\nu\geq 0.$$ We used here that $h^2(\beta)=0$ by Serre duality, combined with the fact that $\beta$ is effective and $p_g=0$. In consequence, for an arbitrary class $\gamma$, we have \begin{eqnarray*} \chi^{\vir}\left(\mathsf {Quot}_{X} (\mathbb C^1, \beta, n), \gamma\right)&=&\chi\left(X^{[m]}\times \mathbb P(V), \,\,\wedge_{-1} \left(M^{[m]}\otimes \zeta^{-1}\right)\otimes \gamma\right)\\&=&\sum_{k=0}^{m} (-1)^k \chi\left(X^{[m]}\times \mathbb P(V), \,\,\wedge^{k} M^{[m]}\otimes \zeta^{-k}\otimes \gamma\right).\end{eqnarray*} Under the identification of the Quot scheme with $X^{[m]}\times \mathbb P$, the tautological classes take the form $$\alpha_{\mathsf {Quot}}^{[n]}=\chi(\alpha)\cdot \mathcal O - \chi(\widetilde \alpha) \cdot \zeta^{-1} + \widetilde{\alpha}_{\mathsf{Hilb}}^{[m]}\otimes \zeta^{-1}$$ where $\widetilde \alpha=\alpha\otimes \mathcal O(-\beta)$; see \cite{JOP}, {Section 3.2}. As a result, we write $$\overline{\mathsf Z}=\mathsf Z_1+\mathsf Z_2+\mathsf Z_3$$ where $$\mathsf Z_1=\chi(\alpha)\cdot \sum_{m=0}^{\infty} q^m \sum_{k=0}^{m} (-1)^k \chi\left(X^{[m]}\times \mathbb P(V), \,\,\wedge^{k} M^{[m]}\otimes \zeta^{-k}\right)$$ $$\mathsf Z_2=-\chi(\widetilde{\alpha})\cdot \sum_{m=0}^{\infty} q^m \sum_{k=0}^{m} (-1)^k \chi\left(X^{[m]}\times \mathbb P(V), \,\,\wedge^{k} M^{[m]}\otimes \zeta^{-k-1}\right)$$ $$\mathsf Z_3=\sum_{m=0}^{\infty} q^m \sum_{k=0}^{m} (-1)^k \chi\left(X^{[m]}\times \mathbb P(V), \,\,\wedge^{k} M^{[m]}\otimes \widetilde{\alpha}^{[m]} \otimes \zeta^{-k-1}\right).$$ 

All three series yield rational functions. The first two series require Theorem 5.2.1 in \cite {Sc} or Corollary 6.1 in \cite{A}: $$\chi\left(X^{[m]}, \,\,\wedge^{k} M^{[m]}\right)=\binom{m-k+\chi(\mathcal O_X)-1}{m-k}\cdot \binom{\chi(M)}{k}=\binom{\nu+1}{k}$$ using that in our case $\chi(\mathcal O_X)=1.$ 
Thus \begin{align*}\mathsf Z_1&=\chi(\alpha)\cdot \sum_{m=0}^{\infty} q^m \sum_{k=0}^{m} (-1)^k \chi\left(X^{[m]}, \,\,\wedge^{k} M^{[m]}\right)\cdot \chi\left(\mathbb P(V), \zeta^{-k}\right)\\&=\chi(\alpha)\cdot \sum_{m=0}^{\infty} q^m \cdot \sum_{k=0}^{m}(-1)^k \binom{\nu+1}{k}\cdot \binom{-k+\nu}{\nu}\\&=\chi(\alpha)\cdot \left(1+q+\ldots+q^\nu\right).\end{align*} The last evaluation follows by examining the nonzero binomial contributions $k=0$ and $k=\nu+1$ in the sum above. 
In a similar fashion, we find $$\mathsf Z_2=-\chi(\widetilde{\alpha})\cdot (\nu+1)\cdot q^{\nu}.$$ Regarding the third sum, we note \begin{align*}\mathsf Z_3&=\sum_{m=0}^{\infty} q^m \sum_{k=0}^{m} (-1)^k \chi\left(X^{[m]}\times \mathbb P(V), \,\,\wedge^{k} M^{[m]}\otimes \widetilde{\alpha}^{[m]} \otimes \zeta^{-k-1}\right)\\&=\sum_{m=0}^{\infty} q^m \sum_{k=0}^{m} (-1)^k \chi\left(X^{[m]}, \wedge^{k} M^{[m]}\otimes \widetilde{\alpha}^{[m]}\right)\cdot \binom{\nu-k-1}{\nu}. \end{align*}
We require the following: 
\begin{proposition}\label{ppp}For all nonsingular surfaces $X$, line bundles $M\to X$, and $K$-theory classes $\widetilde{\alpha}$ on $X$, we have \begin{equation}\label{profor}\mathsf W=\sum_{m=0}^{\infty} q^m \chi\left(X^{[m]}, \wedge_{\,y\,} M^{[m]}\otimes \widetilde{\alpha}^{[m]}\right)=q \cdot \frac{(1+qy)^{\chi(M)}}{(1-q)^{\chi(\mathcal O_X)}} \cdot \chi\left(X, \frac{\wedge_{\,y\,} M}{\wedge_{\,qy\,} M}\otimes \widetilde{\alpha}\right).\end{equation}\end{proposition} \noindent To prove rationality of $\mathsf Z_3$, note that $$\mathsf Z_3=\frac{1}{\nu !}\cdot \frac{\partial^\nu\mathsf W}{\partial y^\nu}\bigg\vert_{y=-1}.$$ We examine the derivatives of each term in the product \eqref{profor}. For instance, for $\nu\geq 1$, we have $$\frac{\partial^\nu}{\partial y^\nu}\chi\left(X, \frac{\wedge_{\,y\,} M}{\wedge_{\,qy\,} M}\otimes \widetilde \alpha\right)\bigg\vert_{y=-1}=(-1)^{\nu-1}\,\nu!\, q^{\nu-1} (1-q) \cdot \chi\left(X, \, \frac{M^{\nu}}{(1-qM)^{\nu+1}}\otimes \widetilde \alpha\right).$$ Expanding using the binomial theorem, we arrive at the expression $$(-1)^{\nu-1}\,\nu!\, q^{\nu-1} (1-q) \cdot \sum_{k=0}^{\infty} (-q)^k \cdot \binom{-\nu-1}{k}\cdot \chi(X, \, M^{\nu+k}\otimes \widetilde \alpha).$$ By Hirzebruch-Riemann-Roch, the Euler characteristics $\chi(X, \, M^{\nu+k}\otimes \widetilde \alpha)$ depend polynomially on $k$, and the same is true about the binomial prefactors. Thus, the answer is a rational function in $q$ with possible pole only at $q=1$. 

The interested reader can compute $\mathsf Z_3$ explicitly by following the above method, substituting $\chi(M)=\chi(\beta)=\nu+1$ and $\chi(\mathcal O_X)=1.$ For example, when $\nu\geq 2$, after cancellations, we obtain $$\mathsf Z_3=q^{\nu}\cdot  (-\chi(M^{-1}\otimes \widetilde \alpha)+ (\nu+1)\cdot \chi(\widetilde \alpha))$$ so that \begin{equation}\label{zbarpg}\overline{\mathsf Z}=\chi(\alpha)\cdot (1+q+\ldots+q^{\nu})-\chi(K^{-1}\otimes \alpha)\cdot q^{\nu}.\end{equation} The answers for $\nu=0$ or $1$ have poles at $q=1$ of order at most $2$;  we do not record them here.\qed 

\vskip.05in
\noindent {\it Proof \footnote{This also follows from \cite {Z}, equation (194) by setting $u=-y,$ $v=0$ in the corresponding formulas. In \cite {Z}, the last term is claimed to be $\chi(X, \widetilde{ \alpha})$. There is however a substitution oversight on the previous page just above equation (192). When the correct values of $u, v$ are used, the above expression naturally appears. For the convenience of the reader, we decided to include an independent argument here.} of Proposition \ref{ppp}.} Since both sides are additive in $\widetilde{\alpha}$, it suffices to assume $\widetilde{\alpha}$ is a line bundle. Consider the more general expression $$\mathsf W^{\dagger}=\sum_{n=0}^{\infty} q^n \chi\left(X^{[n]}, \wedge_{\,y\,} M^{[n]}\otimes \text{Sym}_{\,z\,} \widetilde\alpha^{[n]}\right).$$ We show that $$\mathsf W^{\dagger}=\frac{(1+qy)^{\chi(M)}}{(1-q)^{\chi(\mathcal O_X)}}\cdot \left(1+qz\chi\left(X, \frac{\wedge_{\,y\,} M}{\wedge_{\,qy\,} M}\otimes \widetilde \alpha\right)+\ldots\right),$$ where the omitted terms have order at least $2$ in $z$. The Proposition follows from here by extracting the $z$-coefficient. Note that the $z$-constant term $$\mathsf W^{\dagger}\bigg\vert_{z=0}=\sum_{n=0}^{\infty} q^n \chi\left(X^{[n]}, \wedge_{\,y\,} M^{[n]}\right)=\frac{(1+qy)^{\chi(M)}}{(1-q)^{\chi(\mathcal O_X)}}$$ is indeed correct by \cite{A, Sc}, and was already used above. Consequently, we only need to confirm the asymptotics of the normalized expression \begin{equation}\label{wstar}\mathsf W^{\star}=\frac{\mathsf W^{\dagger}}{\mathsf W^{\dagger}\big\vert_{z=0}}\stackrel{?}{=}1+qz\chi\left(X, \frac{\wedge_{\,y\,} M}{\wedge_{\,qy\,} M}\otimes \widetilde \alpha\right)+\ldots.\end{equation}

By the universality results of \cite {EGL}, it suffices to establish \eqref{wstar} for toric surfaces $X$. Via equivariant localization in $K$-theory, we reduce to the case of equivariant $X=\mathbb C^2$. These steps are standard, and the interested reader can find more details in a similar context in \cite{A, EGL}. 

Let the torus $\mathbb C^{\star}\times \mathbb C^{\star}$ act on $\mathbb C^2$ with weights $t_1, t_2$, which we write multiplicatively. The line bundles $M$ and $\widetilde \alpha$ are trivial, but they carry possibly nontrivial equivariant weights $m$ and $a$ respectively. The equivariant version of \eqref{wstar} becomes $$\mathsf W^{\star}=1+qz\cdot \frac{a}{(1-t_1^{-1})(1-t_2^{-1})}\cdot \frac{1+ym}{1+qym}+\ldots.$$ The torus fixed points of the Hilbert scheme $\left(\mathbb C^2\right)^{[n]}$ are standardly known to correspond to monomial ideals; these are indexed by partitions $\lambda$ of $n$. We note the equivariant restrictions of the tautological bundles to the fixed points $$M^{[n]}\bigg\vert_{\lambda}=\sum_{(b_1, b_2)\in \lambda} m\, t_1^{-b_1}t_2^{-b_2}, \quad \alpha^{[n]}\bigg\vert_{\lambda}=\sum_{(b_1, b_2)\in \lambda} a\, t_1^{-b_1}t_2^{-b_2},$$ as well as the character of the tangent space $$T_{\lambda} \left(\mathbb C^2\right)^{[n]}=\sum_{\square\in \lambda} \left(t_1^{-\ell (\square)}t_2^{a(\square)+1}+t_1^{\ell(\square)+1} t_2^{-a(\square)}\right).$$ The equivariant sum 
$$\mathsf W^{\dagger}=\sum_{n=0}^{\infty} q^n \chi\left(\left(\mathbb C^2\right)^{[n]}, \wedge_{\,y\,} M^{[n]}\otimes \text{Sym}_{\,z\,} \widetilde\alpha^{[n]}\right)$$ thus takes the form $$\mathsf W^{\dagger}=\sum_{\lambda} \frac{q^{|\lambda|}}{\wedge_{-1} T_{\lambda}^{\vee}} \cdot \prod_{(b_1, b_2)\in \lambda}\frac{1+ym\, t_1^{-b_1} t_2^{-b_2}}{1-za\, t_1^{-b_1}t_2^{-b_2}}.$$ It is more convenient to introduce the following modified expression $$\mathsf F(q, u, v)=\sum_{\lambda} \frac{q^{|\lambda|}}{\wedge_{-1} T_{\lambda}^{\vee}} \cdot \prod_{(b_1, b_2)\in \lambda}\frac{1-u t_1^{b_1} t_2^{b_2}}{v- t_1^{b_1}t_2^{b_2}}.$$ Then $$\mathsf W^{\dagger}=\mathsf F\left(-qym, -\frac{1}{ym}, za\right), \quad \mathsf W^{\star}=\frac{\mathsf F\left(-qym, -\frac{1}{ym}, za\right)}{\mathsf F\left(-qym, -\frac{1}{ym}, 0\right)}.$$ 

Crucially for our arguments, $\mathsf F$ admits the following $(q,v)$-symmetry $$\frac{\mathsf F(q, u, v)}{\mathsf F(q, u, 0)}=\frac{\mathsf F(v, u, q)}{\mathsf F(v, u, 0)}.$$ This is contained in Theorem 1.2 of \cite {A}, and is proven via localization over the Donaldson-Thomas moduli space of Calabi-Yau toric threefolds. The symmetry greatly simplifies our calculation: by switching parameters, we only need to find the contribution of the Hilbert scheme of one point. Specifically, we have $$\mathsf W^{\star}=\frac{\mathsf F\left(-qym, -\frac{1}{ym}, za\right)}{\mathsf F\left(-qym, -\frac{1}{ym}, 0\right)}=\frac{\mathsf F\left(za, -\frac{1}{ym}, -qym\right)}{\mathsf F\left(za, -\frac{1}{ym}, 0\right)}.$$ From the definitions, we see that $$\mathsf F\left(za, -\frac{1}{ym}, -qym\right)=1+\frac{za}{(1-t_1^{-1})(1-t_2^{-1})} \cdot \frac{1+\frac{1}{ym}}{-qym-1}+\text{h.o.t in }z,$$ $$\mathsf F\left(za, -\frac{1}{ym}, 0\right)=1-\frac{za}{(1-t_1^{-1})(1-t_2^{-1})}\cdot \left(1+\frac{1}{ym}\right)+\text{h.o.t in }z,$$ from where the equivariant version of \eqref{wstar} follows at once. \qed

\subsection{Rank $1$ quotients} We end this section with a calculation involving rank $1$ quotients. Let $X$ be a nonsingular projective surface with $p_g=0$. Let $\mathsf {Quot}_n$ parametrize short exact sequences \begin{equation}\label{sese}0\to S\to \mathbb C^2\otimes \mathcal O_X\to Q\to 0\end{equation} where $$\text{rk }Q=1,\quad c_1(Q)=0,\quad \chi(Q)=n+\chi(\mathcal O_X).$$ It was shown in \cite{OP, Sch} that $\mathsf {Quot}_n$ admits a virtual fundamental class of dimension $$\text{vd\,}=\chi(\mathcal O_X)$$ whenever $p_g=0$. We write $\chi=\chi(\mathcal O_X)=1-h^1(\mathcal O_X)$ for simplicity. We must have $0\leq \chi\leq 1$ in order to obtain non-zero invariants. As an application of our  computations in Section \ref{srse}, we show
 
\begin{proposition} The series $$\mathsf Z=\sum_{n=0}^{\infty} q^n \chi^{\vir}(\mathsf {Quot}_n, \mathcal O)$$ is the Laurent expansion of a rational function in $q$. 
\end{proposition} 
The method of proof is well-suited to handle twists by tautological sheaves as well. It is natural to inquire whether the Proposition can be generalized to higher $N$ and higher rank quotients $Q$, and to include curve classes $\beta\neq 0$. These questions are more difficult and we will address them elsewhere. 

\proof We let the torus $\mathbb C^\star$ act with weights $w_1, w_2$ on the middle term of the sequence \eqref{sese}. Set $w=w_2-w_1$ for simplicity. There are two fixed point loci $$\mathsf F^\pm=X^{[n]}$$ corresponding to exact sequences $$0\to S=I_Z\to \mathcal O_X[w_1]+\mathcal O_X[w_2]\to Q=\mathcal O_Z+\mathcal O_X\to 0$$ where $\text{length}(Z)=n.$
The $\pm$ decorations correspond to the inclusion of the subsheaf $S$ into the first or second trivial summands respectively. The fixed part of the tangent-obstruction theory is 
$$T^{\vir} \mathsf F^\pm=\text{Ext}^{\bullet}_X(I_{\mathcal Z}, \mathcal O_{\mathcal Z})=TX^{[n]}-\text{Obs}=TX^{[n]}-\left(K^{[n]}\right)^{\vee}.$$ Here $\mathcal Z$ is the universal subscheme. Thus, the fixed locus $\mathsf F^{\pm}=X^{[n]}$ comes equipped with the virtual fundamental class we studied above. Similarly, the virtual normal bundles carry weights $\pm w$ and are given by 
\begin{align*}\mathsf N^\pm&= \text{Ext}_X^{\bullet}(I_{\mathcal Z}, \mathcal O_X)=\text{Ext}_X^{\bullet}(\mathcal O_X-\mathcal O_{\mathcal Z}, \mathcal O_X)=\mathbb C^{\chi} \otimes \mathcal O- \text{Ext}_X^{\bullet}(\mathcal O_{\mathcal Z}, \mathcal O_X)\\&=\mathbb C^{\chi}\otimes \mathcal O-\left(K^{[n]}\right)^{\vee},\end{align*} using Serre duality in the last step. In particular, for $y=e^{-w}$, we have 
$$\frac{1}{\wedge_{-1} \mathsf N^{+\vee}}=\frac{\wedge_{-1} \left(K^{[n]}[-w]\right)}{\left(\wedge_{-1} \mathbb C[-w]\right)^{\chi}}=\left(\frac{1}{1-y}\right)^{\chi}\cdot \left(\sum_{k} (-1)^k \wedge^k K^{[n]} y^k\right)
=\frac{\wedge_{-y\,} K^{[n]}}{(1-y)^{\chi}}$$ and similarly $$\frac{1}{\wedge_{-1} \mathsf N^{-\vee}}=\frac{\wedge_{-y^{-1}\,} K^{[n]}}{(1-y^{-1})^{\chi}}.$$

By virtual localization in $K$-theory \cite{FG, Qu}, we find \begin{align*}\chi^{\vir}(\mathsf {Quot}_n, \mathcal O)&=\chi^{\vir}\left(\mathsf F^{+}, \frac{1}{\wedge_{-1} \mathsf N^{+\vee}}\right)+\chi^{\vir}\left(\mathsf F^{-}, \frac{1}{\wedge_{-1} \mathsf N^{-\vee}}\right)\\&=\frac{1}{(1-y)^{\chi}}\chi^{\vir}\left(X^{[n]}, {\wedge_{-y\,} K^{[n]}}\right)+\frac{1}{(1-y^{-1})^{\chi}} \chi^{\vir}\left(X^{[n]}, \wedge_{-y^{-1}\,} K^{[n]}\right)\bigg\vert_{y=1}.\end{align*} Expressions of this form were studied in Section \ref{srse}. In particular, from the formulas for the functions $\mathsf A$ and $\mathsf B$ immediately following equation \eqref{changev}, we find $$\sum_{n=0}^{\infty} q^n\chi^{\vir}\left(X^{[n]}, \wedge_{-y\,} K^{[n]}\right)=\left(\frac{1-y+yu^2}{1-y+yu}\right)^{K^2}$$ where 
$u$ is the solution of the equation $$-q=\frac{u}{1-u}\cdot \frac{1}{1-y+yu},\quad u(q=0)=0.$$ To evaluate the second term, we similarly solve the equation
 $$-q=\frac{v}{1-v}\cdot \frac{1}{1-y^{-1}+y^{-1}v},\quad v(q=0)=0.$$
 Therefore, 
\begin{align*}
\mathsf Z&=\sum_{n=0}^{\infty} q^n \chi^{\vir}(\mathsf {Quot}_n, \mathcal O)\\&=\frac{1}{(1-y)^{\chi}}\cdot \left(\frac{1-y+yu^2}{1-y+yu}\right)^{K^2}+ \frac{1}{(1-y^{-1})^{\chi}}\cdot \left(\frac{1-y^{-1}+y^{-1}v^2}{1-y^{-1}+y^{-1}v}\right)^{K^2}\bigg\vert_{y=1}.
\end{align*} 

In order to cancel the poles at $y=1$ in the above expression, we record the following:
\begin{lemma}\label{llll}We have 
$$u=(1-y)\cdot \frac{-q}{1+q} + (1-y)^2\cdot \frac{-q^2(2+q)}{(1+q)^3}+O\left((1-y)^3\right)$$
$$v=(1-y)\cdot \frac{q}{1+q} + (1-y)^2 \cdot \frac{q}{(1+q)^3}+ O\left((1-y)^3\right).$$
\end{lemma} 
\proof We only give the details of the first expansion. We could derive this directly by implicit differentiation as in Section \ref{srse}, but a different and slightly stronger argument will be given instead. (The same ideas can be used to give a different proof of Theorem \ref{n1}.) We change variables $1-y=t$, so that $$-q=\frac{u}{1-u}\cdot \frac{1}{t+(1-t)u}.$$ We claim that the implicit solution takes the form \begin{equation}\label{expp}u(q, t)=\sum_{n\geq 1, k\geq 1} a_{n, k}\,q^n t^k \end{equation} where $$a_{n, k}=\frac{(-1)^{n+k+1}}{n} \binom{n+k}{k-1} \binom{n}{k}.$$ The Lemma follows since $$a_{n, 1}=(-1)^n,\quad a_{n, 2}=(-1)^{n+1}\cdot \frac{(n-1)(n+2)}{2},$$ and $$\sum_{n=1}^{\infty} a_{n, 1}\, q^n=\frac{-q}{1+q},\quad \sum_{n=1}^{\infty} a_{n, 2}\, q^n=-\frac{q^2(2+q)}{(1+q)^3}.$$ Expansion \eqref{expp} uses Lagrange inversion for implicit functions \cite{So}, Theorem 3.6. Specifically, let $$\mathsf G(u, q, t)=-q(1-u)(t+(1-t)u)$$ so that we implicitly solve $\mathsf G(u, q, t)=u$. Therefore, $$u(q, t)=\sum_{n\geq 1, k\geq 1} a_{n, k}q^n t^k$$ where \begin{align*}a_{n,k}&= \sum_{m} \frac{1}{m} [z^{m-1}\,q^n\, t^k] \,{\mathsf G}(z, q, t)^m\\&=\sum_{m} \frac{1}{m} [z^{m-1} \,q^n\,t^k]\,(-1)^m \,q^m \,(1-z)^m \,(t+(1-t) z)^m\\&=\frac{(-1)^n}{n}[z^{n-1} \,t^k] \,(1-z)^n \,(t+(1-t)z)^n \\&=\frac{(-1)^n}{n}\, [z^{n-1}\,t^k] \,(1-z)^n (z+t(1-z))^n\\&=\frac{(-1)^n}{n}\,\binom{n}{k}\, [z^{n-1}]\, (1-z)^n\, z^{n-k}\, (1-z)^k\\&=\frac{(-1)^n}{n}\binom{n}{k} [z^{k-1}] (1-z)^{n+k}=\frac{(-1)^{n+k+1}}{n} \binom{n}{k} \binom{n+k}{k-1}.\end{align*}\qed

Using Lemma \ref{llll}, we find that 
$$\frac{u}{1-y}=\frac{-q}{1+q} + (1-y)\cdot \frac{-q^2(2+q)}{(1+q)^3}+O\left((1-y)^2\right),$$ and from here we derive
$$\frac{1-y+yu^2}{1-y+yu}=\frac{1+\frac{yu^2}{1-y}}{1+\frac{yu}{1-y}}=(1+q)\left(1-(1-y)\cdot \frac{q(1-q)}{(1+q)^2}+O\left((1-y)^2\right)\right).$$ A similar argument shows  
$$\frac{1-y^{-1}+y^{-1}v^2}{1-y^{-1}+y^{-1}v}=(1+q)\cdot \left(1+(1-y)\cdot \frac{q(1-q)}{(1+q)^2}+O\left((1-y)^2\right)\right).$$ The Proposition follows from here. For instance, when $\chi=0$, corresponding to surfaces of irregularity $1$, we immediately obtain $$\mathsf Z=2(1+q)^{K^2}.$$ For regular surfaces, we have $\chi=1$. In this case, we need to consider the linear terms in $1-y$ as well, yielding \begin{align*}\mathsf Z&=\frac{(1+q)^{K^2}}{1-y} \left(1-K^2(1-y)\frac{q(1-q)}{(1+q)^2}\right)+\frac{(1+q)^{K^2}}{1-y^{-1}} \left(1+K^2(1-y)\frac{q(1-q)}{(1+q)^2}\right)\bigg\vert_{y=1}\\&=(1+q)^{K^2} \cdot \left(1-2K^2 \cdot \frac{q(1-q)}{(1+q)^2}\right).\end{align*} This completes the proof.\qed

\section{K-theoretic cosection localization}\label{s4}

To complement the explicit calculations of Section \ref{s3}, we now discuss the virtual structure sheaf of Quot schemes with $\beta=0$, giving a proof of Theorem \ref{tt}.  

As in the Introduction, assume $C\subset X$ is a nonsingular canonical curve, and let $\Theta=N_{C/X}=\iota^{\star} K_X$ denote the theta characteristic. Let $$\iota:\mathsf {Quot}_C(\mathbb C^N, n)\hookrightarrow \mathsf {Quot}_X(\mathbb C^N, n)$$ be the natural inclusion. We show that the virtual structure sheaf of $\mathsf {Quot}_X(\mathbb C^N, n)$ localizes on the curve $C$: \begin{equation}\label{col}\mathcal O_{\mathsf{Quot}_X(\mathbb C^N, n)}^{\vir}=(-1)^{n} \,\iota_{\star}\, \mathcal D_n\end{equation}
where we set $$\mathcal D_n\to \mathsf {Quot}_{C}(\mathbb C^N, n),\quad \mathcal D_n =\det {\textnormal{R}}\pi_{\star}(\Theta\otimes \det \mathcal Q)^{\vee}.$$ Along the way, we also prove that \begin{equation}\label{vsq}\mathcal D_n^{\otimes 2} =\det \mathcal N^{\vir}\end{equation} where $\mathcal N^{\vir}$ stands for the virtual normal bundle of the embedding $\iota$. 

To align \eqref{col} with the statement of Theorem \ref{tt} given in the Introduction, consider the support map $$p: \mathsf {Quot}_{C}(\mathbb C^N, n)\to C^{[n]},\,\,\, Q\mapsto \text{supp }Q.$$ 
 We have $$\mathcal D_n=p^{\star}\det \left(\Theta^{[n]}\right).$$ Indeed, for $\mathcal Z\subset C\times C^{[n]}$ denoting the universal subscheme, over $C^{[n]}$ we establish that $$\det {\textnormal R}\pi_{\star} (\Theta\otimes \mathcal O(\mathcal Z))^{\vee}=\det {\textnormal R}\pi_{\star} (\Theta\otimes \mathcal O_{\mathcal Z}).$$ This follows by relative duality, the fact that $\Theta$ is a theta characteristic, and the ideal exact sequence for $\mathcal Z$: $$\det {\textnormal R}\pi_{\star} (\Theta\otimes \mathcal O(\mathcal Z))^{\vee}=\det {\textnormal R}\pi_{\star} (\Theta\otimes \mathcal O(-\mathcal Z))^{\vee}=\det {\textnormal R}\pi_{\star} (\Theta\otimes\mathcal O_{\mathcal Z}).$$
 \begin{example} \label{r1} For $N=1$, 
equation \eqref{col} states
$$\mathcal O_{X^{[n]}}^{\vir}=(-1)^n\iota_{\star} \det \Theta^{[n]},$$ for $\iota:C^{[n]}\hookrightarrow X^{[n]}.$ This is indeed correct since $\Obs=\left((K_X)^{[n]}\right)^{\vee}$ as noted in \cite [Section 4.2.3]{OP}, and we have \begin{eqnarray}\label{eee}\mathcal O_{X^{[n]}}^{\text{vir}}&=&\wedge_{-1}\Obs^{\vee}=(-1)^n\wedge^{n} \Obs^{\vee}\otimes \wedge_{-1} \Obs\\&=&(-1)^n\det (K_X)^{[n]}\otimes \iota_{\star} \mathcal O_{C^{[n]}}\notag\\&=&(-1)^n\iota_{\star} \det \Theta^{[n]}.\notag\end{eqnarray}  Here, we used the Koszul resolution $$\wedge^{\bullet} \left((K_X)^{[n]}\right)^{\vee}\to \mathcal O_{C^{[n]}}$$ of the canonical section $s^{[n]}$ cutting out $\iota:C^{[n]}\hookrightarrow X^{[n]},$ where $s$ defines $\iota:C\hookrightarrow X.$
Furthermore, \begin{eqnarray*}\mathcal N^{\vir} &=&\iota^{\star} T^{\text{vir}}X^{[n]}-TC^{[n]}\\&=&\iota^{\star}\left(TX^{[n]}-\text{Obs}\right)-TC^{[n]}\\
&=&N_{C^{[n]}/X^{[n]}}-\iota^{\star} \text{Obs}\\&=&\Theta^{[n]}-\left(\Theta^{[n]}\right)^{\vee}.\end{eqnarray*} On the last line, we used the identifications $$N_{C^{[n]}/X^{[n]}}=\Theta^{[n]}, \quad  \iota^{\star}\text{Obs}=\iota^{\star} \left((K_C)^{[n]}\right)^{\vee}=\left(\Theta^{[n]}\right)^{\vee},$$ which can be found in \cite [Section 4.3.2]{OP}. This allows us to confirm equation \eqref{vsq}: $$\det \mathcal N^{\vir} = \left(\det \Theta^{[n]}\right)^2.$$ 
\end{example} 
 
\noindent {\it Proof of Theorem \ref{tt}.} We justify equations \eqref{col}  and \eqref{vsq}. For simplicity, we write $\mathsf {Quot}_C$ and $\mathsf {Quot}_X$ for the two Quot schemes. We let $\mathcal S_C, \mathcal S_X$ denote the two universal subsheaves, and let $\mathcal Q$ be universal quotient for $\mathsf {Quot}_C$. We begin by noting that \begin{align*}\mathcal N^{\vir}&=\iota^{\star} T^{\vir} \mathsf {Quot}_X - T\mathsf {Quot}_C\\&=\text{Ext}_X^{\bullet}(\mathcal S_X, \mathcal Q) - \text{Ext}_C^{\bullet} (\mathcal S_C, \mathcal Q)\\&= - \text{Ext}_X^{\bullet}(\mathcal Q, \mathcal Q) + \text{Ext}_C^{\bullet}(\mathcal Q, \mathcal Q) \\&=\text{Ext}^{\bullet}_C (\mathcal Q, \mathcal Q\otimes \Theta).
\end{align*} The third equality follows by expressing the universal subsheaves in terms of the universal quotient in $K$-theory. The last identity follows from the exact sequence $$\ldots \to \text{Ext}_C^i(\mathcal Q, \mathcal Q)\to \text{Ext}_X^i(\mathcal Q, \mathcal Q)\to \text{Ext}_C^{i-1}(\mathcal Q, \mathcal Q\otimes \Theta)\to \ldots $$ proven, for instance, in \cite[Lemma $3.42$]{T}. 

For future reference, over $\mathsf {Quot}_C(\mathbb C^N, n)$, we set $$\mathcal U_n = \det \mathcal N^{\vir}=\det \text{Ext}^{\bullet}_C (\mathcal Q, \mathcal Q\otimes \Theta).$$ We show first that \begin{equation}\label{col1}\mathcal O_{\mathsf{Quot}_X(\mathbb C^N, n)}^{\vir}=(-1)^{n} \,\iota_{\star}\, \mathcal U_n^{1/2}\end{equation} where the right hand side is only well-defined up to $2$-torsion. 

To this end, we apply $\mathbb C^\star$-equivariant localization to both Quot schemes over $C$ and $X$ using the same weights  $w_1, \ldots, w_N$ for the two torus actions. The fixed loci are $$\mathsf F_C[n_1, \ldots, n_N]=C^{[n_1]}\times \cdots \times C^{[n_N]}\, ,
\ \ \  \mathsf F_X[n_1, \ldots, n_N]=X^{[n_1]}\times \cdots \times X^{[n_N]}$$ where $n_1+\ldots+n_N=n$. There is a natural embedding $$\iota: \mathsf F_C[n_1, \ldots, n_N]\hookrightarrow \mathsf F_X[n_1, \ldots, n_N].$$ In the localization argument below, we will match contributions coming from the two fixed loci corresponding to the same values of $n_1, \ldots, n_N$. For simplicity, we drop the $n_i$'s from the notation, writing $\mathsf F_C$ and $\mathsf F_X$ for the fixed loci. We let $$j_C: \mathsf F_C\to \mathsf {Quot}_C, \,\,\,\, j_X: \mathsf F_X\to \mathsf{Quot}_X$$ denote the natural inclusions. We also write $\mathsf N_C$ and $\mathsf N_X$ for the virtual normal bundles. 

The proof of \eqref{col1} requires several steps. By the virtual localization theorem in K-theory \cite{Qu}, we have $$\mathcal O_{\mathsf{Quot}_{X}}^{\mathrm{vir}}=\sum (j_X)_{\star} \,\,\,\, \frac{\mathcal O_{\mathsf{F}_X}^{\vir}}{\wedge_{-1}\mathsf N_X^{\vee}}$$ $$\mathcal O_{\mathsf{Quot}_{C}}=\sum (j_C)_{\star}\,\,\, \frac{\mathcal O_{\mathsf{F}_C}}{\wedge_{-1}\mathsf N_C^{\vee}}.$$ In particular, from the second expression, we find  
\begin{align*}(-1)^n\iota_{\star}\, \mathcal U_n^{1/2} &= \sum (-1)^{n} \iota_{\star} (j_C)_{\star}\,\,\, \frac{j_C^{\star}\,\, \mathcal U_n^{1/2}}{\wedge_{-1}\mathsf N_C^{\vee}}\\&=\sum (-1)^{n}  (j_X)_{\star}\,\iota_{\star}\,\,\, \frac{j_C^{\star}\,\, \mathcal U_n^{1/2}}{\wedge_{-1}\mathsf N_C^{\vee}}.\end{align*} To prove \eqref{col1}, it suffices to establish\begin{equation}
\label{sp} \frac{\mathcal O_{\mathsf {F}_X}^{\vir}}{\wedge_{-1}\mathsf N_X^{\vee}}=(-1)^n \iota_{\star}\,\,\, \frac{j_C^{\star}\,\, \mathcal U_n^{1/2}}{\wedge_{-1}\mathsf N_C^{\vee}}.\end{equation} 

We compute the virtual structure sheaf on the left. The obstruction theory of the fixed locus $\mathsf {\mathsf F}_X=X^{[n_1]}\times \cdots \times X^{[n_N]}$ splits as $$\Obs=\oplus_{i=1}^{N} \Obs_i,\quad \Obs_i=\left(K_X^{[n_i]}\right)^{\vee}.$$ Thus, using \eqref{eee}, we find \begin{equation}\label{e2}\mathcal O_{\mathsf{F}_X}^{\vir}=\wedge_{-1} \Obs^{\vee} = \otimes_{i=1}^{N} \wedge_{-1}\Obs_i^{\vee} =\iota_{\star}\left( \otimes_{i=1}^{N} (-1)^{n_i} \det \Theta^{[n_i]}\right).\end{equation} To prove \eqref{sp}, we show \begin{equation}\label{e1} \frac{\otimes_{i=1}^{N} \det \Theta^{[n_i]}}{\iota^{\star}\wedge_{-1}  \mathsf N_X^{\vee}} = \frac{j_C^{\star}\,\, \mathcal U_n^{1/2}}{\wedge_{-1}\mathsf N_C^{\vee}}.\end{equation}

Next, we note that over the fixed locus $\mathsf F_C$, the quotient takes the form $$\mathcal Q=\mathcal O_{\mathcal Z_1}+\ldots+\mathcal O_{\mathcal Z_N}.$$ Thus the equivariant restriction becomes 
\begin{equation}\label{e3}j_{C\,}^{\star} \mathcal U_n=\det  \text{Ext}_{C}^{\bullet}\left(\sum_{i=1}^{N} \mathcal O_{\mathcal Z_i}, \sum_{i=1}^{N} \mathcal O_{\mathcal Z_i}\otimes \Theta\right).\end{equation} Both sides carry weight $0$. 

Lemma 34 in \cite {OP} gives the difference \begin{equation}\label{diff}\iota^{\star} \mathsf N_X-\mathsf N_C=\bigoplus_{i<j} \mathsf V_{ij}\end{equation}  where we define $$\mathsf V_{ij}=\text{Ext}^{\bullet}_C( \mathcal O_{\mathcal Z_i}, \mathcal O_{\mathcal Z_j}\otimes \Theta)[w_j-w_i] \oplus \text{Ext}^{\bullet}_C( \mathcal O_{\mathcal Z_j}, \mathcal O_{\mathcal Z_i}\otimes \Theta)[w_i-w_j].$$  We noted in \cite {OP} that by Serre duality and the fact that $\Theta$ is a theta characteristic, we have $$\mathsf V_{ij}=\mathsf U_{ij}+\mathsf U_{ij}^{\vee}[-1]\, $$ where $$\mathsf U_{ij}=\text{Ext}^{\bullet}_C( \mathcal O_{\mathcal Z_i}, \mathcal O_{\mathcal Z_j}\otimes \Theta)[w_j-w_i] .$$ Observe that $$\text{rank }\mathsf U_{ij}=0.$$

In general, if $\mathsf V=\mathsf U+ \mathsf U^{\vee}[-1]$ we have \begin{equation} \label{e4}\frac{1}{\wedge_{-1} \mathsf V^{\vee}}= (-1)^{\text{rk } \mathsf U}\det \mathsf U.\end{equation} 
To see this, note that since both sides are multiplicative, we may assume $\mathsf U$ is an equivariant line bundle $L$. We then need to show $$\frac{1}{\wedge_{-1} \left(L^{\vee} - L\right)}=-L\iff \wedge_{-1}L=-\wedge_{-1} L^{\vee} \cdot L$$ which is correct. 

From \eqref{diff} and \eqref{e4} we find 
\begin{equation}\label{e5}\frac{1} {\iota^{\star}\wedge_{-1}\mathsf N_X^{\vee}}=\frac{1}{\wedge_{-1} N_C^{\vee}} \cdot \prod_{i<j} \frac{1}{\wedge_{-1} \mathsf V_{ij}^{\vee}}=\frac{1}{\wedge_{-1} N_C^{\vee}} \cdot \det \left(\sum_{i<j} \mathsf U_{ij}\right).\,\end{equation} Comparing \eqref{e3} and \eqref{e5} with \eqref{e1}, we need to show \begin{equation}\label{ddddd}\otimes_{i} \det \Theta^{[n_i]}\otimes \det \left(\sum_{i<j} \mathsf U_{ij}\right)=\det  \left(\text{Ext}_{C}^{\bullet}\left(\sum_{i=1}^{N} \mathcal O_{\mathcal Z_i}, \sum_{i=1}^{N} \mathcal O_{\mathcal Z_i}\otimes \Theta\right)\right)^{1/2}.\end{equation} All terms above carry zero weights. A priori, the bundles $\mathsf U_{ij}$ carry weights $w_j-w_i$, but the ranks are equal to $0$, so the weights cancel after taking determinants. 

To establish \eqref{ddddd} we extend the definition of $\mathsf U_{ij}$ to include the case $i=j$: $$\mathsf U_{ii}=\text{Ext}^{\bullet}_C( \mathcal O_{\mathcal Z_i}, \mathcal O_{\mathcal Z_i}\otimes \Theta).$$ By Serre duality, we have  \begin{align*}\det \mathsf U_{ii}&=\det\left(\text{Ext}_C^0(\mathcal{O}_{\mathcal{Z}_i}, \mathcal O_{\mathcal Z_i}\otimes \Theta)-\text{Ext}_C^1(\mathcal{O}_{\mathcal{Z}_i}, \mathcal O_{\mathcal Z_i}\otimes \Theta)\right)\\&=\det\left(\text{Ext}_C^0(\mathcal{O}_{\mathcal{Z}_i}, \mathcal O_{{\mathcal Z}_i}\otimes \Theta)-\text{Ext}_C^0(\mathcal{O}_{\mathcal{Z}_i}, \mathcal O_{\mathcal Z_i}\otimes \Theta)^{\vee}\right)\\&=\det \text{Ext}_C^0(\mathcal{O}_{\mathcal{Z}_i}, \mathcal O_{{\mathcal Z}_i}\otimes \Theta)^{\otimes 2}=\det \text{Ext}_C^0(\mathcal O_C, \mathcal O_{\mathcal Z_i}\otimes \Theta)^{\otimes 2}\\&=\det {\Theta^{[n_i]}}^{\otimes 2}.\end{align*} The ideal exact sequence for $\mathcal Z_i$ was used on the third line. 
We already noted above that as a consequence of Serre duality $\mathsf U_{ij}\simeq \mathsf U_{ji}^{\vee}[1]$, hence $$\det \mathsf U_{ij}=\det \mathsf U_{ji}.$$ Therefore \begin{align*} \otimes_i \det \Theta^{[n_i]}\otimes \det  \left(\sum_{i < j} \mathsf U_{ij}\right)&=\det\left(\sum_{1\leq i,j \leq N} \mathsf U_{ij}\right)^{1/2}\\&=\det\left(\text{Ext}^{\bullet}_C\left(\sum_{i=1}^{N} \mathcal{O}_{\mathcal{Z}_i},\sum_{i=1}^{N} \mathcal{O}_{\mathcal{Z}_i}\otimes \Theta\right )\right)^{1/2}.\end{align*}This confirms \eqref{ddddd}, and equation \eqref{col1} along with it. 

Finally, recall that $$\mathcal D_n, \,\,\mathcal U_n\to \mathsf {Quot}_{C}(\mathbb C^N, n),\quad \mathcal D_n =\det {\textnormal{R}}\pi_{\star}(\det \mathcal Q\otimes \Theta)^{\vee}, \quad  \mathcal U_n=\det \text{Ext}_C^{\bullet}(\mathcal Q, \mathcal Q\otimes \Theta).$$ We show that $$\mathcal D_n^{\otimes 2} =\mathcal U_n.$$ If in the expression for $\mathcal U_n$ we replace $\mathcal Q$ by $\det \mathcal Q-\mathcal O,$ we obtain \begin{align*}\widetilde{\mathcal U}_n &= \det \text{Ext}^{\bullet}_C (\det \mathcal Q- \mathcal O, \det \mathcal Q\otimes \Theta- \Theta)\\&=\det \left(\textnormal{R}\pi_{\star}(\det \mathcal Q\otimes \Theta)\right)^{\vee} \otimes \det \textnormal{R}\pi_{\star}(\det \mathcal Q^{\vee}\otimes \Theta)^{\vee}\\&=\det \left(\textnormal{R}\pi_{\star}(\det \mathcal Q\otimes \Theta)\right)^{\vee\otimes 2}=\mathcal D_n^{\otimes 2},\end{align*} where Serre duality was used in the last line. It remains to show that $\mathcal U_n\simeq \mathcal {\widetilde U}_n$. Write $$\Delta= (\det \mathcal Q-\mathcal O)-\mathcal Q,$$ and note that $$\mathcal {\widetilde U}_n= \mathcal U_n\otimes \det \text{Ext}_C^{\bullet}(\Delta, \Delta\otimes \Theta)\otimes \det \text{Ext}_C^{\bullet}(\Delta, \mathcal Q\otimes \Theta)\otimes \det \text{Ext}_C^{\bullet}(\mathcal Q, \Delta\otimes \Theta).$$ With respect to the codimension filtration in $K$-theory, $\Delta$ has pieces of degree $2$ or higher, and $\mathcal Q$ has pieces of degree $1$ or higher. It follows that the last three Ext's above are supported in codimension $2$, hence their determinants are trivial. This completes the argument. \qed

\section{The Segre/Verlinde correspondence and symmetry}\label{s5}
\subsection{Overview}
We study here the virtual shifted Verlinde and Segre series \begin{eqnarray*} \mathsf V_{\alpha} (q)&=&q^{\beta. K_X}\sum_{n=0}^{\infty} q^n \, \chi^{\textrm{vir}}(\mathsf{Quot}_{X}(\mathbb C^N, \beta, n),\, \det \alpha^{[n]})\\\mathsf S_{\alpha} (q)&=&q^{\beta. K_X}\sum_{n=0}^{\infty} q^n \int_{\left[\mathsf{Quot}_{X}(\mathbb C^N, \beta, n)\right]^{\textrm{vir}}} s(\alpha^{[n]}).\end{eqnarray*}
We establish the Segre/Verlinde correspondence of Theorem \ref{t2}. For $\beta=0$, we prove the symmetry of the Segre/Verlinde series with respect to $N$ and $r=\text{rank }\alpha$, as stated in Theorem \ref{ss}. 
 \subsection{The Verlinde Series} \label{vser} The virtual Verlinde series corresponds to the series $\mathsf Z_{X, N, \beta}^{f, g}(\alpha)$ of Section \ref{s1} for 
$$f(x)=e^x, \quad g(x)=\frac{x}{1-e^{-x}}.$$ This is a consequence of the virtual Hirzebruch-Riemann-Roch theorem \cite {CFK, FG}. 
In particular, 
$$\frac{f'}{f}(x)=1, \quad \frac{g(x)}{x}=\frac{1}{1-e^{-x}}, \quad \frac{g'}{g}(-x)+\frac{1}{x}=\frac{1}{1-e^{-x}}.$$
We can simultaneously simplify the universal series and the change of variables by setting 
\begin{equation}\label{ex}e^{H_i}=1-\widetilde{H}_i \text{ and }
e^{-w_k}=1+\widetilde{w}_k,\end{equation}
so that equation \eqref{chvar} becomes 
\[q=\big(1-\widetilde{H}_i\big)^r\cdot \prod_{k=1}^{N} \frac{\widetilde{H}_i+\widetilde{w}_k}{1+\widetilde{w}_k}\quad \text{with}\quad \widetilde{H}_i(q=0)=-\widetilde{w}_i.
\]
Using \eqref{bv} -- \eqref{W}, we obtain the explicit expressions 
\begin{align*}
\mathsf{A}&={
{\prod\limits_{i=1}^{N} \left(\frac{1-\widetilde{H}_i}{1+\widetilde{w}_i}\right)^{-r}\cdot \prod\limits_{k, i} \frac{\widetilde{H}_i+\widetilde{w}_k}{-(1-\widetilde{H}_i)}\cdot \prod\limits_{i_1\neq i_2} \frac{1-\widetilde{H}_{i_1}}{\widetilde{H}_{i_2}-\widetilde{H}_{i_1}}}\cdot \prod\limits_{i=1}^{N} \left(r-\sum\limits_{k=1}^{N} \frac{1-\widetilde{H}_i}{\widetilde{H}_i+\widetilde{w}_k}\right)}\\
\mathsf{B}&=\prod_{i=1}^{N} \frac{1-\widetilde{H}_i}{1+\widetilde{w}_i}\\
\mathsf{U}_i&=
\left(1-\widetilde{H}_i\right)^r\cdot \prod\limits_{k=1}^N \frac{-(1-\widetilde{H}_i)}{\widetilde{H}_i+\widetilde{w}_k}\cdot \prod\limits_{i'\neq i} \left(\frac{\widetilde{H}_{i'}-\widetilde{H}_{i}}{1-\widetilde{H}_{i}}\cdot  \frac{\widetilde{H}_{i}-\widetilde{H}_{i'}}{1-\widetilde{H}_{i'}}\right)\cdot {\left(r-\sum\limits_{k=1}^{N} \frac{1-\widetilde{H}_i}{\widetilde{H}_i+\widetilde{w}_k}\right)^{-1}}\\
\mathsf{V}_i&=\frac{1}{1-\widetilde{H}_i}\\
\mathsf{W}_{i_1, i_2}&=\frac{1-\widetilde{H}_{i_1}}{\widetilde{H}_{i_2}-\widetilde{H}_{i_1}}\cdot  \frac{1-\widetilde{H}_{i_2}}{\widetilde{H}_{i_1}-\widetilde{H}_{i_2}}.
\end{align*}

\subsection{The Segre Series}\label{sser} To obtain the virtual Segre series, we consider the functions 
\[\widetilde {f}(x)=\frac{1}{1+x}, \quad \widetilde {g}(x)=1.
\]
In particular, 
\[\frac{\widetilde {f}'}{\widetilde {f}}(x)=-\frac{1}{1+x},  \quad \frac{\widetilde {g}}{x}=\frac{1}{x}, \quad \frac{\widetilde {g}'}{\widetilde{g}}(-x)=0.
\]
We write $\widetilde w_1, \ldots, \widetilde w_N$ for the weights used in the corresponding localization computation. Via equations \eqref{chvar} -- \eqref{W}, we find the universal series
\begin{align*}
\widetilde{\mathsf{A}}&= \prod\limits_{i=1}^N\left(\frac{1-\widetilde{H}_i}{1+\widetilde{w}_i}\right)^{-r}\cdot \prod\limits_{k, i}(\widetilde H_i+\widetilde w_k)\cdot \prod\limits_{i_1\neq i_2}\frac{1}{\widetilde H_{i_1}-\widetilde H_{i_2}} \cdot {\prod\limits_{i=1}^{N} \left(\frac{-r}{1-\widetilde{H}_i}+\sum\limits_{k=1}^{N} \frac{1}{\widetilde{H}_i+\widetilde{w}_k}\right)}\\
\widetilde{\mathsf{B}}&=\prod_{i=1}^{N}\frac{1-\widetilde{H}_i}{1+\widetilde{w}_i}\\
\widetilde{\mathsf{U}}_i&={\left(1-\widetilde{H}_i\right)^{r}\cdot \prod\limits_{k=1}^N\frac{1}{\widetilde H_i+\widetilde w_k}\cdot \prod\limits_{i'\neq i}\big(\widetilde H_{i'}-\widetilde H_{i}\big)\big(\widetilde H_{i}-\widetilde H_{i'}\big) }\cdot {\left(\frac{-r}{1-\widetilde{H}_i}+\sum\limits_k \frac{1}{\widetilde{H}_i+\widetilde{w}_k}\right)^{-1}}\\
\widetilde{\mathsf{V}}_i&=\frac{1}{1-\widetilde{H}_i}\\
\widetilde{\mathsf{W}}_{i_1, i_2}&=\frac{1}{\widetilde H_{i_1}-\widetilde H_{i_2}}\cdot \frac{1}{\widetilde H_{i_2}-\widetilde H_{i_1}}\end{align*}
after the change of variables 
\[\widetilde{q}=\big(1-\widetilde{H}_i)^r \cdot \prod_{k=1}^{N} \big(-\widetilde{H}_i-\widetilde{w}_k\big) \ \ \text{with} \ \ \widetilde{H}_i(q=0)=-\widetilde{w}_i.
\]

\subsection{Comparison} We show the Segre and Verlinde series match\footnote{There are other quadruples $(f, g, \widetilde f, \widetilde g)$ for which the corresponding series of invariants match as well, but their geometric interpretation is less clear.} for the following geometries
\begin{itemize}
\item [(i)] $X$ is smooth projective, $\beta=0$, $N$ is arbitrary, 
\item [(ii)] $X$ is smooth projective with $p_g>0$, $N=1$ and $\beta$ is arbitrary, 
\item [(iii)] $X$ is a relatively minimal elliptic surface, $\beta$ is supported on fibers, $N$ is arbitrary.
\end{itemize} 

\noindent {\it Proof of Theorem \ref{t2}.} We begin by rewriting the Segre series as $${\mathsf {S}}^{\dagger}_{\alpha} (q)=q^{\beta. K_X}\sum_{n=0}^{\infty} q^n \int_{\left[\mathsf{Quot}_{X}(\mathbb C^N, \beta, n)\right]^{\textrm{vir}}} s\left(\alpha^{[n]}-\left(\mathbb C^{N}\right)^{\oplus n}\right).$$ This does not change the original nonequivariant expression. Equivariantly, we use the lift of the trivial bundle $\mathbb C^N$ with weights $\widetilde w_1, \dots, \widetilde w_N$.
Assuming that the torus weights $w_1, \ldots, w_N$ on the Verlinde side  and the torus weights $\widetilde w_1, \dots, \widetilde w_N$ on the Segre side are connected by \eqref{ex}, we show that equivariantly $$\mathsf V_{\alpha}(q)={\mathsf S}^{\dagger}_{\alpha}((-1)^{N} q).$$ Theorem \ref{t2} follows from here. 

We begin by analyzing the series $\mathsf V_{\alpha}$ and $\mathsf S_{\alpha}$. By the calculations in Sections \ref{vser} and \ref{sser}, the following universal series match exactly on both series:
\[\mathsf{A}=\widetilde{\mathsf{A}}, \quad \mathsf{B}=\widetilde{\mathsf{B}}, \quad \mathsf{V}_i=\widetilde{\mathsf{V}}_i
\]
when evaluated at the variables $q$ and $\tilde q$ related by \begin{equation}\label{change}\frac{q}{\widetilde{q}}=(-1)^N \cdot \prod_{k=1}^{N} \frac{1}{1+\widetilde{w}_k}.
\end{equation} 
However, the remaining universal series no longer match precisely. In fact, 
\begin{align*}
\frac{\mathsf{U}_i}{\widetilde{\mathsf{U}}_i}&=(-1)^{N-1} \cdot \prod_{i'\neq i} \frac{1}{1-\widetilde{H}_{i'}}\\
\frac{\mathsf{W}_{i_1,i_2}}{\widetilde{\mathsf{W}}_{i_1,i_2}}&=(1-\widetilde{H}_{i_1})(1-\widetilde{H}_{i_2}).
\end{align*}
The geometries (i)--(iii) prevent these terms from taking effect: in cases (i) and (iii) they appear with exponent $0$, while in case (ii) only $\mathsf U$ makes sense and in this case we do obtain a perfect match. 
Thus, for (i)--(iii), we obtain via \eqref{usus} the equivariant identity 
$$\mathsf V_{\alpha}(q)=\mathsf S_{\alpha}(\widetilde {q})$$ where $$\frac{q}{\widetilde{q}}=(-1)^N \cdot \prod_{k=1}^{N} \frac{1}{1+\widetilde{w}_k}.$$
We can absorb the difference in the variables by working with ${\mathsf {S}}^{\dagger}_{\alpha}$ instead. Recall the lift of the trivial bundle $\mathbb C^N$ with weights $\widetilde w_1, \dots, \widetilde w_N$. This affects each individual summand in ${\mathsf S}^{\dagger}_{\alpha}$ by the factor $$(1+\widetilde {w}_1)^{n}\cdots (1+\widetilde{w}_N)^{n}\,\,\text { so that }\,\,\mathsf S_{\alpha}\left(q\cdot \prod_{k=1}^{N} ({1+\widetilde{w}_k})\right)={\mathsf S}^{\dagger}_{\alpha}(q).$$ The claimed equality follows $$\mathsf V_{\alpha}(q)={\mathsf S}^{\dagger}_{\alpha}((-1)^{N} q).$$\qed
\begin{example} When $N=1$ and for $\text{rank }\alpha=r$, we have 
\begin{equation}\label{ee}\sum_{n=0}^{\infty} q^n \cdot \chi^{\textrm{vir}}(X^{[n]},\, \det \alpha^{[n]})=(1-t)^{c_1(\alpha)\cdot K_X} \cdot \left(\frac{1-t(r+1)}{(1-t)^{r+1}}\right)^{K_X^2}\end{equation} for $q=t(1-t)^r.$ This is a consequence of Theorem \ref{t2} and \cite[Theorem 2]{MOP2} giving the Segre series. Thus, the theorems in Section \ref{kthry} do not extend to the highest exterior powers and to the Verlinde series which are typically given by algebraic functions. 

Rationality does occur however in the special case when $\text{rank }\alpha=0$ (also for $\text{rank }\alpha=-1$), so that $\det \alpha^{[n]}$ agrees with the pullback from the symmetric power $\left(\det \alpha\right)_{(n)}.$
Then, 
\begin{equation}\label{eef}\sum_{n=0}^{\infty} q^n \cdot \chi^{\textrm{vir}}(X^{[n]},\, \det \alpha^{[n]})=(1-q)^{-\chi^{\vir}(X, \alpha)}=(1-q)^{c_1(\alpha)\cdot K_X}\, .\end{equation} 
By \eqref{eee}, we have $$\mathcal O_{X^{[n]}}^{\text{vir}}=\sum_{k=0}^{n}(-1)^k \wedge^k (K_X)^{[n]}\, .$$
The computation \eqref{eef} is thus aligned with Theorem 5.2.1 in \cite {Sc} which gives the individual cohomology groups $$H^\star(X^{[n]}, (\det \alpha)_{(n)} \otimes \wedge^k (K_X)^{[n]})=\wedge^kH^{\star}(\det \alpha\otimes K_X)\otimes \mathsf{Sym}^{n-k} H^{\star}(\det \alpha)\, .$$ 
\end{example}

\subsection{Symmetry} We now prove Theorem \ref{ss} giving the symmetry of the Segre/Verlinde series when $\beta=0$. We pick two classes $\alpha$ and $\widetilde \alpha$ such that $$\text{rk }\alpha=r,\quad \text{rk }\widetilde \alpha=N,\quad\frac{c_1(\alpha)\cdot K}{\text{rk }\alpha}=\frac{c_1(\widetilde{\alpha})\cdot K}{\text{rk }{\widetilde \alpha}}=\mu.$$ We show $$\mathsf S_{N, \,\alpha}((-1)^{N}q)=\mathsf S_{r, \,\widetilde{\alpha}}((-1)^{r}q).$$

For clarity, we change the notation from the previous subsections, writing $$\mathsf S_{N, \,\alpha}((-1)^{N}q)=\sum_{n=0}^{\infty} \left((-1)^{N}q\right)^n\int_{\left[\mathsf{Quot}_{X}(\mathbb C^N, n)\right]^{\vir}} s(\alpha^{[n]})=\mathsf M^{K^2} \cdot \mathsf N^{\mu}$$ where \begin{align*}\mathsf M&={\prod\limits_{i=1}^{N} \left(\frac{-r}{1-{H}_i}+\sum\limits_{k=1}^{N} \frac{1}{H_i+{w}_k}\right)}\, {\prod\limits_{i=1}^N\left(\frac{1-{H}_i}{1+{w}_i}\right)^{-r} \prod\limits_{1\leq i,\, k\leq N}{(H_i+ w_k)}\prod\limits_{i_1\neq i_2}\frac{1}{H_{i_1}-H_{i_2}}}\\
\mathsf{N}&=\prod_{i=1}^{N} \left(\frac{1-{H}_i}{1+{w}_i}\right)^{r}\end{align*} for \begin{equation}\label{bam}q=\big(1-{H}_i)^r \cdot \prod_{k=1}^{N} \big({H}_i+{w}_k\big),\quad H_i(q=0)=-w_i.\end{equation} In a similar fashion, $$\mathsf S_{r,\,\widetilde \alpha}((-1)^{r}q)=\sum_{n=0}^{\infty} \left((-1)^{r}q\right)^n\int_{\left[\mathsf{Quot}_{X}(\mathbb C^r, n)\right]^{\vir}} s(\widetilde \alpha^{[n]})=\widetilde{\mathsf M}^{K^2} \cdot \widetilde{\mathsf N}^{\mu}$$ where  \begin{align*}\widetilde{\mathsf M}&={\prod\limits_{j=1}^{r} \left(\frac{-N}{1-\widetilde {H}_j}+\sum\limits_{k=1}^{r} \frac{1}{\widetilde H_j+\widetilde {w}_k}\right)}\, {\prod\limits_{j=1}^r\left(\frac{1-\widetilde {H}_j}{1+\widetilde {w}_j}\right)^{-N} \prod\limits_{1\leq j, \,k\leq r}{(\widetilde H_j+ \widetilde w_k)}\prod\limits_{j_1\neq j_2}\frac{1}{\widetilde H_{j_1}-\widetilde H_{j_2}}}\\
\widetilde{\mathsf{N}}&=\prod_{i=1}^{r} \left(\frac{1-\widetilde {H}_j}{1+\widetilde {w}_j}\right)^{N}\end{align*} for $$q=\big(1-\widetilde {H}_j)^N \cdot \prod_{k=1}^{r} \big(\widetilde {H}_j+\widetilde {w}_k\big),\quad \widetilde H_j(q=0)=-\widetilde w_j.$$ Two sets of weights $w=(w_1,\ldots, w_N)$ and $\widetilde w=(\widetilde w_1, \ldots, \widetilde w_r)$ are used on the two different Quot schemes in the localization computation. We will show that $$ \mathsf M\bigg\rvert_{w=0}=\widetilde {\mathsf M}\bigg\rvert_{\widetilde w=0},\quad \mathsf N\bigg\rvert_{w=0}=\widetilde {\mathsf N}\bigg\rvert_{\widetilde w=0}.$$
In fact, we argue that the limit is found by setting $w=0$ directly in the above equations. In other words $$ \mathsf M\bigg\rvert_{w=0}=\mathsf M^{\star},\quad \mathsf N\bigg\rvert_{w=0}=\mathsf N^{\star}$$ where
\begin{align*}\mathsf M^{\star}&={\prod\limits_{i=1}^{N} \left(\frac{-r}{1-H_i}+ \frac{N}{H_i}\right)}\cdot {\prod\limits_{i=1}^N\left({1-H_i}\right)^{-r}\cdot \prod\limits_{i=1}^{N}H_i^N\cdot \prod\limits_{i_1\neq i_2}\frac{1}{H_{i_1}-H_{i_2}}}\\
\mathsf N^{\star}&=\prod_{i=1}^{N} (1-H_i)^r,\end{align*} and $H_1, \ldots, H_N$ are the roots\footnote{By slight abuse, we do not introduce new notation for the roots of the equation obtained by setting $w=0$.} of the equation $$q=(1-H)^r \cdot H^{N}$$ that vanish at $q=0$. These roots are given by formal Puiseux series in the variable $q^{\frac{1}{N}}.$ 
To justify the limit, note that expressions $\mathsf M$ and $\mathsf N$ are symmetric in the $H$'s, so they can be recast in terms of the elementary symmetric functions. The assertion follows from the following:\vskip.1in
\noindent {\bf Claim:} The elementary symmetric functions of the solutions $H_1, \ldots, H_N$ of \eqref{bam} are continuous as functions of $w$ at $w=0$. \vskip.1in

\noindent We will justify this shortly. In a similar fashion, in the limit, $\widetilde{\mathsf M}$ and $\widetilde{\mathsf N}$ become 
\begin{align*}\widetilde{\mathsf M}^{\star}&={\prod\limits_{j=1}^{r} \left(\frac{-N}{1-\widetilde {H}_j}+ \frac{r}{\widetilde {H}_j}\right)}\cdot {\prod\limits_{j=1}^r\left({1-\widetilde {H}_j}\right)^{-N}\cdot \prod\limits_{j=1}^{r}\widetilde H_j^r\cdot \prod\limits_{j_1\neq j_2}\frac{1}{\widetilde H_{j_1}-\widetilde H_{j_2}}}\\
\widetilde{\mathsf N}^{\star}&=\prod_{j=1}^{r} (1-\widetilde H_j)^N\end{align*} where $\widetilde H_1, \ldots, \widetilde H_r$ are the roots of $q=(1-\widetilde H)^N\cdot \widetilde H^r$ that vanish at $q=0$. 

We show $$\mathsf M^{\star}=\widetilde{\mathsf M}^{\star},\quad \mathsf N^{\star}=\widetilde{\mathsf N}^{\star}.$$
Let $$f(z)=z^N(1-z)^r-q.$$ We note that $f(z)=0$ has $r+N$ roots namely $$H_1, \ldots, H_N, 1-\widetilde H_1, \ldots, 1-\widetilde H_r.$$
Thus $$f(z)=\prod_{i=1}^{N}(z-H_i) \cdot \prod_{j=1}^{r} (1-z-\widetilde H_j).$$  Using the equations, we can rewrite $$\widetilde{\mathsf N}^{\star}=\prod_{j=1}^{r}\frac{q}{\widetilde H_j^r}\implies \frac{\mathsf N^{\star}}{\widetilde {\mathsf N}^{\star}}=q^{-r} \prod_{i=1}^{N} (1-H_i)^{r}\cdot \prod_{j=1}^{r} \widetilde H_j^r= q^{-r} f(1)^r (-1)^r=1\implies \mathsf N^*=\widetilde{\mathsf N}^*.$$
The computation for the remaining series is similar. We first note $$f'(z)=z^{N} (1-z)^r \cdot \left(\frac{-r}{1-z}+\frac{N}{z} \right).$$ From here, we find \begin{align*}\mathsf M^{\star} {{\mathsf N}^{\star}}^2&={\prod\limits_{i=1}^{N} \left(\frac{-r}{1-H_i}+ \frac{N}{H_i}\right)}\cdot {\prod\limits_{i=1}^N\left({1-H_i}\right)^{r}\cdot \prod\limits_{i=1}^{N}H_i^N\cdot \prod\limits_{i_1\neq i_2}\frac{1}{H_{i_1}-H_{i_2}}}\\&=\prod_{i=1}^{N} f'(H_i)  \cdot \prod_{i_1\neq i_2} \frac{1}{H_{i_1}-H_{i_2}} \\&= \prod_{1\leq i\leq N, 1\leq j\leq r} (1-H_i-\widetilde H_j)= \widetilde{\mathsf M}^{\star}\widetilde{\mathsf N}^{{\star}^2}.\end{align*} This shows $\mathsf M^{\star}=\widetilde{\mathsf M}^{\star}$, completing the proof of Theorem \ref{ss}. \vskip.1in

\noindent {\it Proof of the Claim:} The proof is a modification of the argument of Section $2.2$ of \cite{JOP} which does not apply directly here. We work over the ring $A=\mathbb C(w)[[q]]$ of formal power series whose coefficients are rational functions in $w$. 
For simplicity, set $$\mathsf P(z)=(1-z)^r, \quad \mathsf Q(z)=\prod_{k=1}^{N} (z+w_k),$$ and write $$\mathsf P(z)=(-1)^r\left(z^r+p_1z^{r-1}+\ldots+p_r\right), \quad \mathsf Q(z)=z^N+q_1 z^{N-1}+\ldots+q_N.$$ By assumption, $H_1, \ldots, H_N\in A$ are the power series solutions of $$\mathsf P(H)\cdot \mathsf Q(H)-q=0,\quad H_i(q=0)=-w_i.$$ Let $e_i\in A$ be the $i^{\text{th}}$ elementary symmetric function in $H_1, \ldots, H_N$ times $(-1)^{i}$, and note that from the initial conditions, we have at $q=0$: $$e_i(0)=q_i.$$ We factor in $A[z]$: \begin{align}\label{factor}\mathsf P(z)\cdot \mathsf Q(z)-q&=(-1)^r\cdot \prod_{k=1}^{N}(z-H_k)\cdot \left(z^r+f
_1z^{r-1}+\ldots+f_r\right)\\&=(-1)^r\cdot (z^N+e_1 z^{N-1}+\ldots+e_N)\cdot \left(z^r+f
_1z^{r-1}+\ldots+f_r\right)\notag.\end{align} Setting $q=0$, we find $$\mathsf P(z)\cdot \mathsf Q(z)=(-1)^r \cdot \prod_{k=1} (z-H_k(0)) \cdot \left(z^r+f
_1(0)z^{r-1}+\ldots+f_r(0)\right).$$ From the initial conditions $H_k(0)=-w_k$, we obtain $$f_j(0)=p_j.$$ We use the convention $e_0=f_0=1$. 

We show by induction on $m$ that the coefficients of $q^m$ in $e_i$ and $f_j$ are rational functions in $w$ whose denominators are powers of $$\Delta=\prod_{k=1}^{N} \mathsf P(-w_k).$$ Since $\mathsf P(0)\neq 0$, the substitution $w=0$ is therefore allowed, completing the proof of the {\bf Claim.} 

The base case $m=0$ is a consequence of the above remarks. For the inductive step, write $$e_i=q_i + q^{m+1} e_i^{(m+1)} + \text {other terms},\quad f_j=p_j+ q^{m+1} f_j^{(m+1)}+\text{ other terms}.$$ The omitted terms correspond to powers of $q$ with exponent between $1$ and $m$, or greater than $m+1$. From \eqref{factor}, we see that $$[q^{m+1}]\sum_{i+j=k} e_i f_j= [q^{m+1}z^{N+r-k}]\, (-1)^r\left(\mathsf P(z)\mathsf Q(z)-q\right)=0 \text{ or }(-1)^{r+1}.$$ On the other hand, by direct computation $$[q^{m+1}]\sum_{i+j=k} e_i f_j=\sum_{i+j=k} p_j e_i^{(m+1)}+ q_i f_j^{(m+1)}+\text{ other terms}.$$ By the inductive hypothesis, we obtain that $$\sum_{i+j=k} p_j e_i^{(m+1)}+ q_i f_j^{(m+1)}=c_k$$ where $c_k$ is a rational function in $w$ whose denominator is a power of $\Delta$ for $1\leq k\leq N+r$. We regard the above equations as a linear system in $$e_1^{(m+1)}, \ldots, e_N^{(m+1)}, f_1^{(m+1)}, \ldots, f_r^{(m+1)}$$ which we solve in terms of the $c_k$'s by Cramer's rule. To complete the inductive step, we show that the determinant of the matrix of coefficients equals $\Delta$. Indeed, the matrix of coefficients has entries equal to $p_j$ or $q_i$, and zero elsewhere. This is in fact the Sylvester matrix of the two polynomials $(-1)^r\mathsf P$, $\mathsf Q$, and its determinant is therefore the resultant $$\Delta=\text{Res }((-1)^r\mathsf P, \mathsf Q)=\prod_{k=1}^{N} \mathsf P(-w_k).$$

\qed
\begin{example} \label{ex7} The match between expressions (i)-(iv) and (i)$'$-(iv)$'$ in Section \ref{sss} is established by combining Theorems \ref{t2}, \ref{ss}, \ref{tt}. Here, we offer an analogy with the projective space.
Indeed, the following identity can be seen as the most basic instance of the Segre/Verlinde correspondence \begin{equation} \label{sv} (-1)^k\int_{\mathbb P^{k}} s(\mathcal O_{\mathbb P^{k}}(1))^{r+1}=\chi(\mathbb P^{k}, \mathcal O_{\mathbb P^k}(r)),\end{equation} while the simplest rank-section symmetry, when formulated on the Verlinde side, is $$\chi(\mathbb P^k, \mathcal O(r))=\chi(\mathbb P^r, \mathcal O(k)).$$ 

In fact the above statements correspond precisely to the case $n=1,\,\,\, \beta=0, \,\,N \text{ arbitrary},$ so that our results can be interpreted as possible extensions to arbitrary $n$. Note that $$\mathsf {Quot}_{X}(\mathbb C^N, 1) = X \times \mathbb P^{N-1}$$  
in such a fashion that $(x, \zeta)\in X\times \mathbb P^{N-1}$ corresponds to the quotient $\mathbb C^{N}\stackrel{\zeta}{\to}\mathbb C\to \mathbb C_{x}\to 0.$ Since $$\mathcal Q=\mathcal O_{\Delta}\boxtimes \mathcal O_{\mathbb P^{N-1}}(1)$$ over $X\times X\times \mathbb P^{N-1},$ we compute 
$$\text{Ext}^{\bullet}(\mathcal S, \mathcal Q)=TX+T\mathbb P^{N-1}-K_X^{\vee}\implies \Obs=K_X^{\vee}$$ and $$\alpha^{[1]} = \alpha\, \boxtimes\, \mathcal O_{\mathbb P^{N-1}}(1)\implies \det \alpha^{[1]}=\det \alpha \,\boxtimes\, \mathcal O_{\mathbb P^{N-1}}(r).$$ Thus the Verlinde number is $$\chi^{\vir}(\mathsf {Quot}_{X}(\mathbb C^N, 1) , \det \alpha^{[1]})=\chi^{\vir}(X, \det \alpha)\cdot \chi(\mathbb P^{N-1}, \mathcal O_{\mathbb P^{N-1}}(r)).$$ On the Segre side, we find \begin{align*}\int_{\left[\mathsf {Quot}_{X}(\mathbb C^N, 1)\right]^{\vir}} s(\alpha^{[1]})&=\int_{X\times \mathbb P^{N-1}}\mathsf e(K_X^{\vee})\cdot \,s\left(\alpha\boxtimes \mathcal O_{\mathbb P^{N-1}}(1)\right)\\
&=-\left(\int_{X} \mathsf e(K_X^{\vee}) \cdot c_1(\alpha)\right)\cdot \int_{\mathbb P^{N-1}} s(\mathcal O_{\mathbb P^{N-1}}(1))^{r+1}\\&=-\chi^{\vir}(X, \det \alpha)\cdot   \int_{\mathbb P^{N-1}} s(\mathcal O_{\mathbb P^{N-1}}(1))^{r+1}.\end{align*} Equation \eqref{sv} confirms Theorem \ref{t2} in this case. We invite the reader to check that all vertices of the cube in Section \ref{sss} match precisely. 
  \end{example}

\end{document}